\documentclass[journal]{IEEEtran}
\ifCLASSINFOpdf
   \usepackage[pdftex]{graphicx}
\else
\fi
\usepackage{multicol}
\usepackage{multirow}

\usepackage{amsmath}
\usepackage{mathtools}
\usepackage{eurosym}
\DeclareRobustCommand{\officialeuro}{%
  \ifmmode\expandafter\text\fi
  {\fontencoding{U}\fontfamily{eurosym}\selectfont e}}

\usepackage[makeroom]{cancel}

\ifCLASSOPTIONcompsoc
  \usepackage[caption=false,font=normalsize,labelfont=sf,textfont=sf]{subfig}
\else
  \usepackage[caption=false,font=footnotesize]{subfig}
\fi
\usepackage{url}
\begin{document}
\title{Clustering Methods Assessment for Investment in Zero Emission Neighborhoods' Energy System}

%\author{Dimitri~Pinel,
%        Magnus~Korp\r{a}s,~\IEEEmembership{Member ?,~IEEE}
%        and~Karen~B.~Lindberg,~\IEEEmembership{Member ?,~IEEE}}% <-this % stops a space
%\thanks{D. Pinel and M. Korp\r{a}s are with the Department
%of Electrical Engineering, NTNU, Trondheim,
%Norway, e-mail: dimitri.q.a.pinel@ntnu.no}% <-this % stops a space
%\thanks{K. B. Lindberg is with SINTEF Byggforsk.}% <-this % stops a space

\author{Dimitri~Pinel \\
Department of Electrical Power Engineering\\
NTNU, Trondheim, Norway\\
E-mail: dimitri.q.a.pinel@ntnu.no\thanks{Address: Elektrobygget, O. S. Bragstads plass 2E, E, 3rd floor, 7034 Trondheim, Norway}}% <-this % stops a space
\maketitle

\begin{abstract}

 This paper investigates the use of clustering in the context of designing the energy system of Zero Emission Neighborhoods (ZEN). ZENs are neighborhoods who aim to have net zero emissions during their lifetime. While previous work has used and studied clustering for designing the energy system of neighborhoods, no article dealt with neighborhoods such as ZEN, which have high requirements for the solar irradiance time series, include a $CO_2$ factor time series and have a zero emission balance limiting the possibilities. 
 To this end several methods are used and their results compared. The results are on the one hand the performances of the clustering itself and on the other hand, the performances of each method in the optimization model where the data is used. Various aspects related to the clustering methods are tested. The different aspects studied are: the goal (clustering to obtain days or hours), the algorithm (k-means or k-medoids), the normalization method (based on the standard deviation or range of values) and the use of heuristic.
 The results highlight that k-means offers better results than k-medoids and that k-means was systematically underestimating the objective value while k-medoids was constantly overestimating it. When the choice between clustering days and hours is possible, it appears that clustering days offers the best precision and solving time. The choice depends on the formulation used for the optimization model and the need to model seasonal storage.
 The choice of the normalization method has the least impact, but the range of values method show some advantages in terms of solving time.
 When a good representation of the solar irradiance time series is needed, a higher number of days or using hours is necessary. The choice depends on what solving time is acceptable.
\end{abstract}


% \fi
\IEEEpeerreviewmaketitle

\section{Introduction}

For optimization models, complexity and solving time are important elements. Some models require to be solved in a limited amount of time. Optimization models used for the control of processes or the unit commitment problem are examples of this. For other models, it can also simply be linked to the practicality of using the model. 
In order to obtain a solution in the desired amount of time, different approaches have been used. One possibility is to simply add a time limit to the model so that if this time limit is reached, the optimization is stopped. The obtained solution is sub-optimal and the distance from optimality will depend on the model and the time limit.
Reducing the complexity of the model, for example by reducing the number of binary variables or linearizing quadratic constraints, is also a possibility.
Another possibility is to reduce the dimensionality of the problem. Clustering is a common way to achieve this reduction. It allows, from a dataset of any dimensionality, to gather points that are close together into a number of clusters. There exist many possible metrics and methods to perform clustering, and their performance will vary depending on the application.
All of those methods will however give sub-optimal solutions to the original problem. The choice of a method should be made according to the desired precision and run time.

In this paper, we investigate the use of clustering in a mixed integer linear program (MILP) called ZENIT. The goal is to identify which technique performs best for this application regarding the time necessary to solve the model, the optimality gap, and the representation of some time series of particular importance.

ZENIT (Zero Emission Neighborhood (ZEN) Investment Tool) is a program based on optimization that helps design the energy system of neighborhoods in a cost-optimal way and with a goal of having achieved net zero emissions of $CO_2$ in the neighborhood's lifetime. 
It is developed as a part of the research center on Zero Emission Neighborhoods in Smart Cities in Norway. The goal of this center is to research solutions to reduce the emission of neighborhoods in various fields such as architecture, urban planning and materials.

In this paper the focus of the clustering is a reduction of the time dimensionality, i.e. using less timesteps. The dimension of the dataset to cluster depends on the length of the time series used and the number of buildings in the neighborhood. 

\section{State of the art and contribution}

Clustering algorithms have been studied extensively since the 1930's \cite{harold32} and improved since then. The principle of those algorithms is to gather similar observations of a dataset into clusters based on a given metric. The outputs of such algorithms are a list of all original data points and the cluster to which they belong as well as a representative vector for each cluster. Many algorithms exist but, in this paper the focus is on the k-means and the k-medoids algorithms because they are the most commonly use for such applications. Those algorithms differ in the way the representative vector of each cluster are chosen. The k-means algorithm uses a centroid as the representative vector, i.e. the vector with the smallest squared distance to every member of the cluster \cite{macqueen1967}. The k-medoids algorithm chooses the representatives of the clusters by choosing the vector in the original data with the smallest distance to every other members of the cluster \cite{kaufmann87}.

Clustering has been extensively studied for multiple applications in various fields, including power systems.
It has been for example used in the context of grid expansion planning in \cite{HARTEL17}, national energy system planning \cite{PFENNINGER20171}\cite{NAHMMACHER16} and unit commitment models \cite{wogrin16}. 

%papers on clustering in power systems in general, grid investment etc
%List of articles:
%\cite{PFENNINGER20171}
%\cite{wogrin16} (modelling energy storage in a unit commitment model, %kmean on hours)
%\cite{HARTEL17}grid expension planning (systematic sampling, k-means, %k-medoids, hierarchical clustering (ward's linkage), moment matching, 
%\cite{NAHMMACHER16} Clustering as a way to represent long term %time series

Many clustering techniques exist and \cite{PFENNINGER20171} suggests that the best choice depends on the data to process and the model in which they are going to be used. It is thus important to compare different methods in order to find the best choice for our particular needs. It also gives insights in the choice of the number of clusters to use.
Several articles compare, with different approaches, the possible clustering techniques. Among them, \cite{PFENNINGER20171} compares the performance of downsampling, k-means and hierarchical clustering as well as different heuristics and combinations of previously mentioned methods. It finds that for their energy system planning model and in the context of pluri-annual time series, some heuristics appear promising. The clustering is performed on days, with 4 different time series and multiple locations giving a rather large number of dimensions.

For a grid expansion planning problem, \cite{HARTEL17} compares systematic sampling, k-means, k-medoids, hierarchical clustering with Ward's linkage and moment matching. It clusters on hours and 5 dimensions. In this case, hierarchical and k-medoids appear to perform equally well.

Closer to the ZENIT model needs, \cite{SCHUTZ18} compared clustering algorithms (k-means, k-centers, k-medoids, k-medians, monthly averaged days, and seasonal days) to find representative days for a model investing and operating the energy system of a building. It finds k-medoids as the best suited method for this application.

Reference \cite{KOTZUR18} also compares different techniques in the context of different local energy systems (averaging, k-means, k-medoids, hierarchical) for obtaining representative days, 3-days or weeks. It finds that medoids perform better than centroids but recommends overall the use of hierarchical clustering due to the reproducibility of the results.
%Comparative studies if any ?
%List of articles:
%\cite{PFENNINGER20171}
%\cite{SCHUTZ18}
%\cite{HARTEL17}
%\cite{KOTZUR18} comparison of different methods and application to %energy system design

It is also interesting to look at the choices made for other models similar to ZENIT, i.e. model for investment and operation in the energy system of buildings or neighborhoods. Those choices are naturally dependant on factors such as the scale of the neighborhood, the level of detail of the model, the target run time, the machine used to solve the model or its goal: investment and/or operation and in some cases grid layout, but it remains a good indication nonetheless.

Many authors choose to use season based clustering (SBC), where they choose or average the time series to form one representative day for each season \cite{CAPUDER14} or only for the summer, the winter and the mid-season \cite{YANG15}\cite{yokoyama02}\cite{Weber11}. They also have varying choices in terms of number of periods for the chosen days: from hourly (i.e. 24 periods)\cite{yokoyama02}, to twelve \cite{YANG15}\cite{yokoyama02}, or six periods\cite{yokoyama02}\cite{Weber11}.
Similarly, some choose to use one average day per month \cite{harb16}\cite{MORVAJ16}\cite{SCHUTZ17}, or several days per month, such as \cite{MASHAYEKH17} with a week day, a week end day and a peak day per month or \cite{PIACENTINO13} with 2 days of 12 periods each per month.

The exact method used to determine the days is not always clear; \cite{ORTIGA11} points this out and suggests a graphical method using the load duration curves. Another method relying on k-means clustering is proposed in \cite{FAZLOLLAHI14}.

Reference \cite{LI17} uses weekly downsampling to allow the model to run faster and checks the scheduling with a 24h rolling horizon model with hourly resolution.
Complete years with hourly resolution are also used in some models \cite{ASHOURI13}.

Other studies rely on clustering \cite{SCHUTZ18}. Reference \cite{GABRIELLI2018408} suggest a way to keep seasonal storage operation while using design days found with k-means clustering. Similarly \cite{STADLER16} relies on k-medoids clustering to find design days. However, only outside temperature and global irradiance are used, assuming that the other time series are correlated to either of those 2. The other time series are reconstructed from the clusters after the clustering.
K-means is also used in \cite{FLEISCHHACKER19}, where two models are coupled, for providing representative weeks and for providing representative hours. The hours clustering is preceded by the removal of peaks from the time series and followed by their re-introduction.

In this paper different methods of clustering, normalizing and treating peaks are compared in the specific case of ZENIT. In addition, design days and representative hours are compared to find the strengths and weaknesses of each approach. This study stands out from other comparative studies by limiting the number of algorithm used but also considering the choices for normalizing and handling peaks. The Zero Emission context also brings specific problems to overcome. For example, the zero emission balance constraint in the optimization model limits the way one can reduce the number of timesteps. Another example is the strong requirements on the solar irradiance time series due to the importance of PV in the results.

\section{Reduction of the number of timesteps}

%Why we cannot use certain methods (downsampling, or certain heuristics) in our model
%Clustering methodologies kmeans kmedoid, heuristic or not, 
%Data preparation for the clustering especially design days
%Normalisation of the data

Many possibilities exist in order to reduce the number of timesteps in the optimization. However some are not adequate for the model. Downsampling for instance is not well suited. With the downsampling method, the time series are reduced by averaging the values on a certain period of time. A six hours downsampling would average the values of the time series on intervals of six hours, dividing by six the total number of timesteps. This method reduces the precision of the data and is not well suited for applications with renewable energies, which vary rapidly. 
The use of heuristic is often considered, and there are different approaches depending on the application. The heuristic could be reducing the time series to a collection of extreme events found in the time series, such as the hours with the maximum load or the lowest temperature or any combination of such criteria. In the case of ZENIT, this is not an acceptable solution on its own. Despite the reduction of the level of details induced, which could be somewhat overcome by tuning the heuristic chosen, the biggest reason that contraindicates its use is the Zero Emission balance constraint. Indeed, using this constraint requires to take into account every hour in the year, which is difficult with heuristics.
On the contrary, clustering allows the use of the Zero Emission balance. In clustering, an algorithm is used to gather similar timesteps into clusters. Each original timestep is then represented by a cluster. We choose this approach over downsampling and heuristic in order to keep the original time granularity and the use of the emission constraint.

Several clustering algorithm exists and we limit this study to k-means and k-medoids clustering. In addition we consider the use of heuristics in combination to the clustering. This approach is recommended in this kind of optimization application because the clustering alone would likely 'dilute' the extreme events' timesteps, such as the hour with the maximum load, into a cluster represented by a lower value, which would lead to an under-dimensionned solution. A simple heuristic in addition to the clustering allows to correct this. In this paper, the heuristic chosen is the time (day or hour) with the highest total load, defined as the sum of the domestic hot water load (DHW), space heating (SH) and electric load, and the time with the lowest irradiance.
In addition, normalizing the data before clustering can be beneficial \cite{scipy}. Several ways to normalize the data before the clustering algorithm exists and we also consider two options: a normalization based on the range of each time series (\ref{norm_range}) and one based on the standard deviation (\ref{norm_std}).

\noindent
\begin{minipage}[c]{0.6\linewidth}
    \begin{equation} \label{norm_range}
        X'=\frac{X-min(X)}{max(X)-min(X)}
    \end{equation}
\end{minipage} % no space if you would like to put them side by side
\begin{minipage}[c]{0.39\linewidth}
    \begin{equation} \label{norm_std}
        X'=\frac{X}{std(X)}
    \end{equation}
\end{minipage}

Lastly, as mentioned in the literature review, mainly two approaches exist for clustering, one clusters directly the hours, the other focuses on design days. The design day approach uses clustering for selecting representative days in the year and then use the hourly values for each representative day. 
This approach is often favored when storages are modelled. Indeed, because the relation between timesteps inside a day are kept, it allows for daily operation of storages contrary to hours clustering.

\begin{figure}
    \centering
    \includegraphics[width=0.47\textwidth,trim=4 4 4 4,clip]{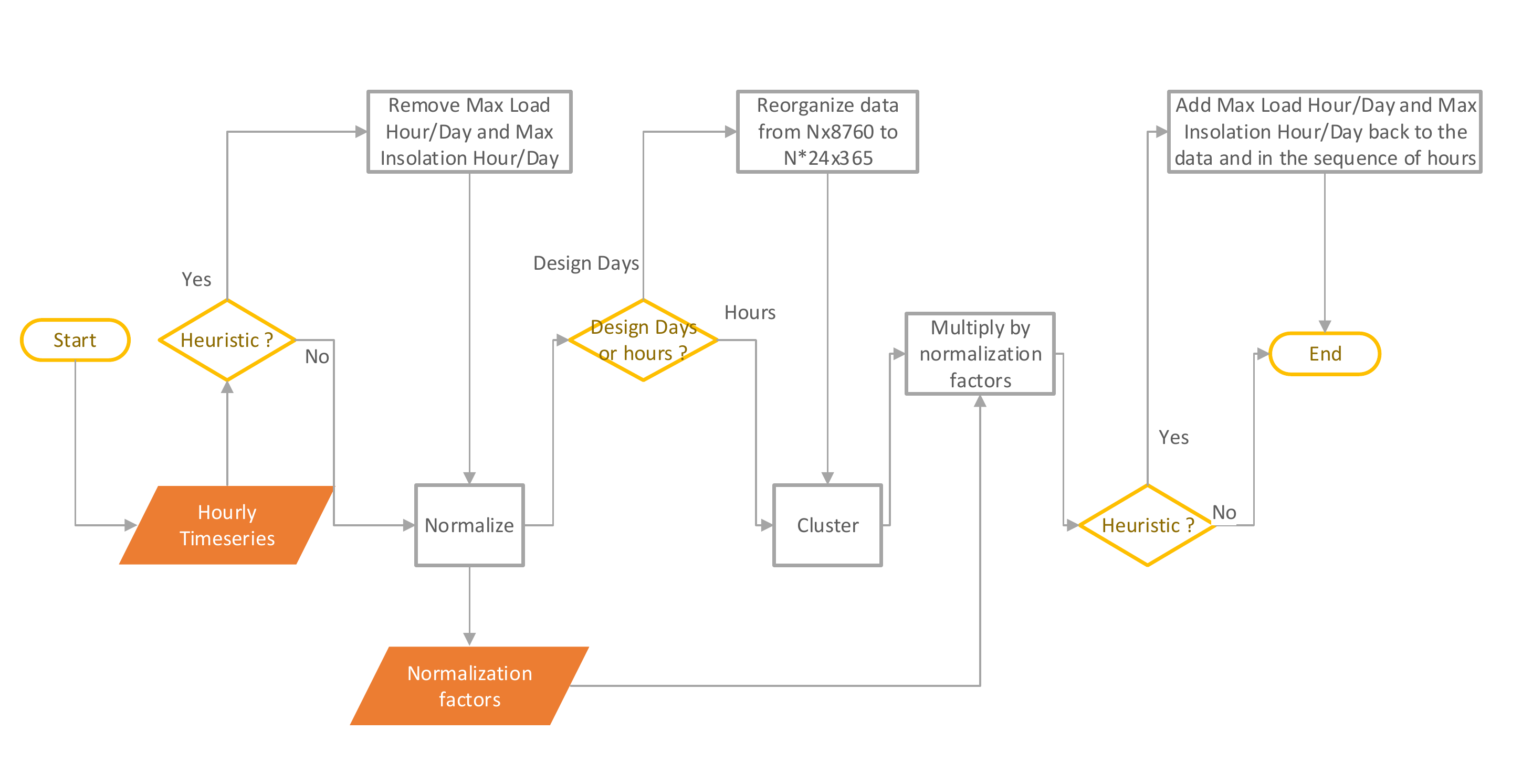}
    \caption{Flowchart of the Clustering Process}
    \label{flowchart_clustering}
\end{figure}

The clustering is performed in Python using PyClustering \cite{Novikov2019} for the k-medoids algorithm and Scipy for the k-means \cite{scipy}\cite{k-means++07}.
The practical handling of the clustering is described in the flowchart in Fig. \ref{flowchart_clustering}.

The data entering the clustering process consists of several hourly time series covering one year. The data is composed of the following time series: one domestic hot water load (DHW), one space heating load (SH) and one electric load for each building (or building type) in the neighborhood; outside air temperature, total irradiance and $CO_2$ factor of electricity.

\section{Clustering results} \label{clust}
\begin{figure}[h]
    \centering
    \includegraphics[width=0.3\textwidth]{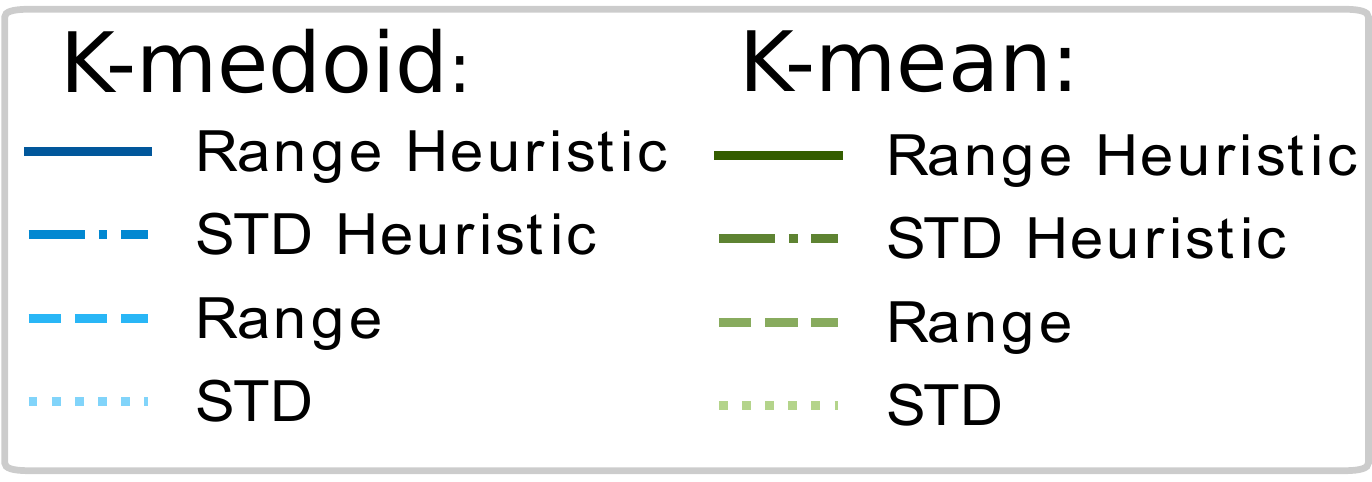}
    \caption{Legend of the Results}
    \label{legend}
\end{figure}

The different clustering approaches presented in the previous section were performed for various number of clusters: for the clustering of design days, up to 100, and for the clustering of hours, up to 2400  (with 6 hours steps). This allows to determine which number of clusters to use in the optimization model.
The representatives of clusters and their sequence are combined to rebuild a complete year and then compared to the original data to compute errors. In this section, the errors are presented as Root Mean Square Deviations (RMSD) and as Normalized RMSD (NRMSD) when comparing the errors of different time series. All figures below share the same legend presented on Fig. \ref{legend}.

\begin{figure}[h]
    \centering
    \subfloat[Hours]{
    \includegraphics[width=0.24\textwidth]{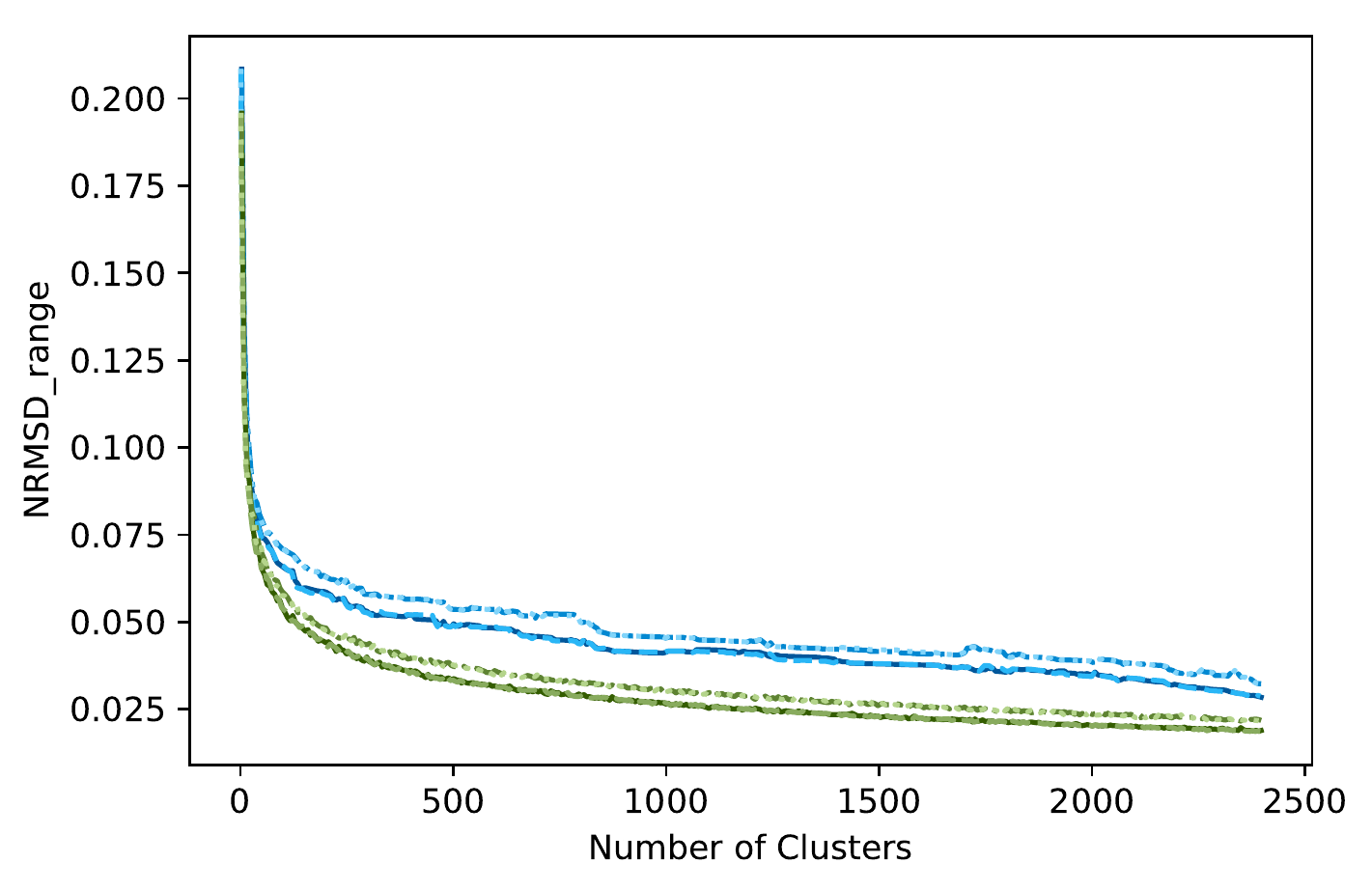}}
    \subfloat[Days]{\includegraphics[width=0.24\textwidth]{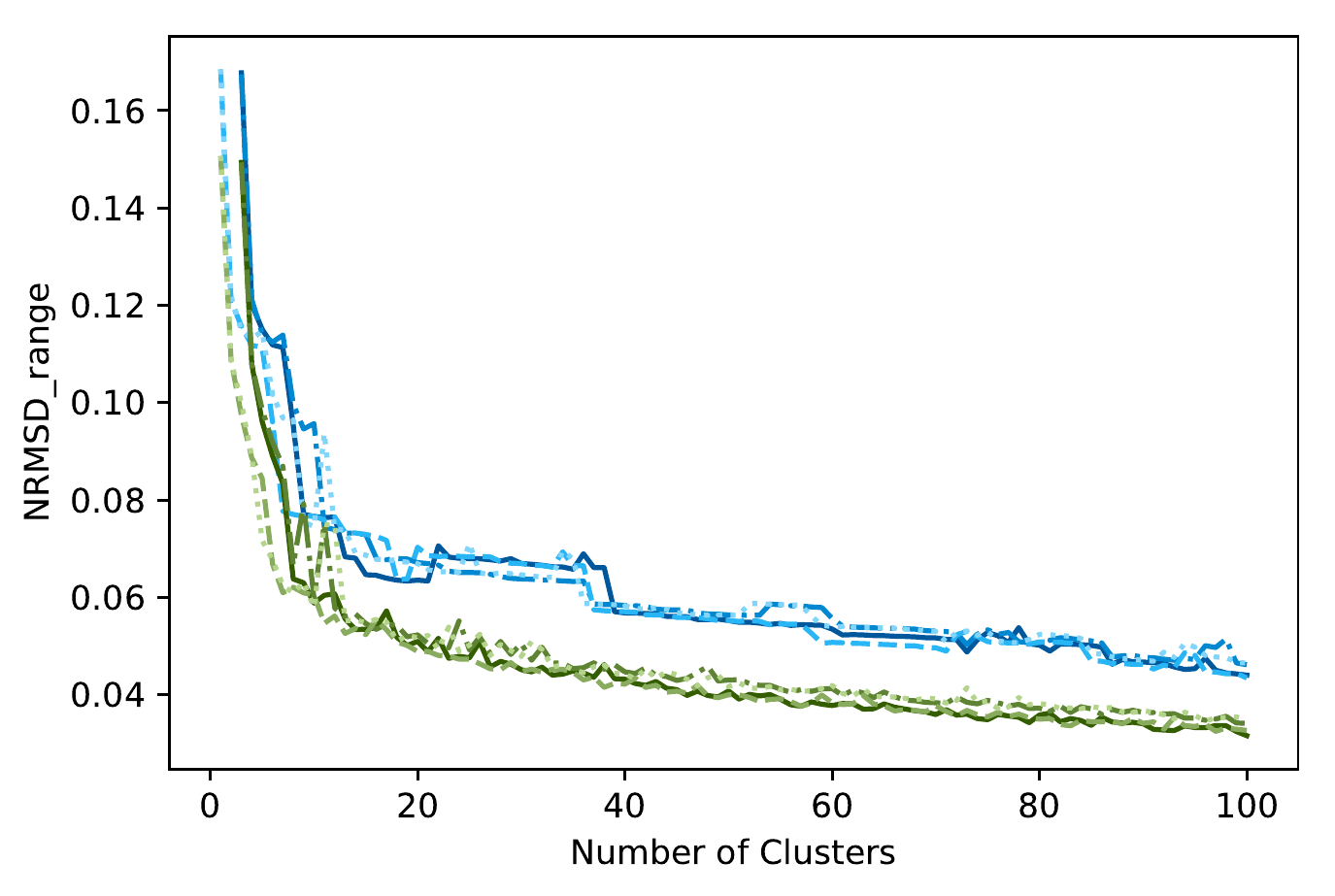}}
    \caption{Average of the NRMSD of All Clustered Time Series, Normalized with Range}
    \label{NRMSD}
\end{figure}

%\begin{figure}[h]
%    \centering
%    \includegraphics[width=0.48\textwidth]{Mean_NRMSD_range_days.pdf}
%    \caption{Average of the NRMSD of all clustered time series, normalized with range, Days}
%\end{figure}
%
%\begin{figure}[h]
%    \centering
%    \includegraphics[width=0.48\textwidth]{Mean_NRMSD_std_hours.pdf}
%    \caption{Average of the NRMSD of all clustered time series, normalized with std, Hours}
%\end{figure}
%
%\begin{figure}[h]
%    \centering
%    \includegraphics[width=0.48\textwidth]{Mean_NRMSD_std_days.pdf}
%    \caption{Average of the NRMSD of all clustered time series, normalized with std, Days}
%\end{figure}

Considering all figures in this section, it is clear that the k-means algorithm performs better than its k-medoids counterparts. This is what we could expect. Indeed, the k-medoids uses vectors from the original datasets instead of creating centroids, which are better representatives. However, this ensures that the chosen representatives of clusters in the case of k-medoids are meaningful and realistic.

Another thing one could expect is that the performance of the clusterings monotonically improves. However, this is not the case of our results, especially in the case of design days. For the performance regarding individual time series, this could be explained because of a better performance of other time series for this particular number of clusters. However, this lack of monotony can still be found in the aggregated result of Fig. \ref{NRMSD}. One possible explanation for the lack of monotony could be that k-means and k-medoids algorithms do not always find the global optimums but can provide solutions that are only local optimums. Hierarchical clustering or running the clustering algorithms several time with different initial conditions could provide more consistent results. 

Looking at Fig. \ref{NRMSD}, the use of heuristic results in a tiny advantage for the heuristic versions on the overall error of the clustering. This is especially true in the case of clustering on hours. For design days clustering, the difference between clustering with and without heuristic disappears after around 8 design days. The lower the amount of design days, the higher the impact of forcing two days to be extreme events is, while for hours, the forced hours are "diluted" faster. 

From all figures, considering an equivalent resulting number of timesteps (translating to the complexity to solve the model) clustering on hours gives much better results than clustering on days. For the overall error, Fig. \ref{NRMSD}, the error for the hours clustering is about 50\% lower than for the design days.

The performance for individual time series is discussed in the following.

\begin{figure}[h]
    \centering
    \subfloat[Hours]{\includegraphics[width=0.24\textwidth]{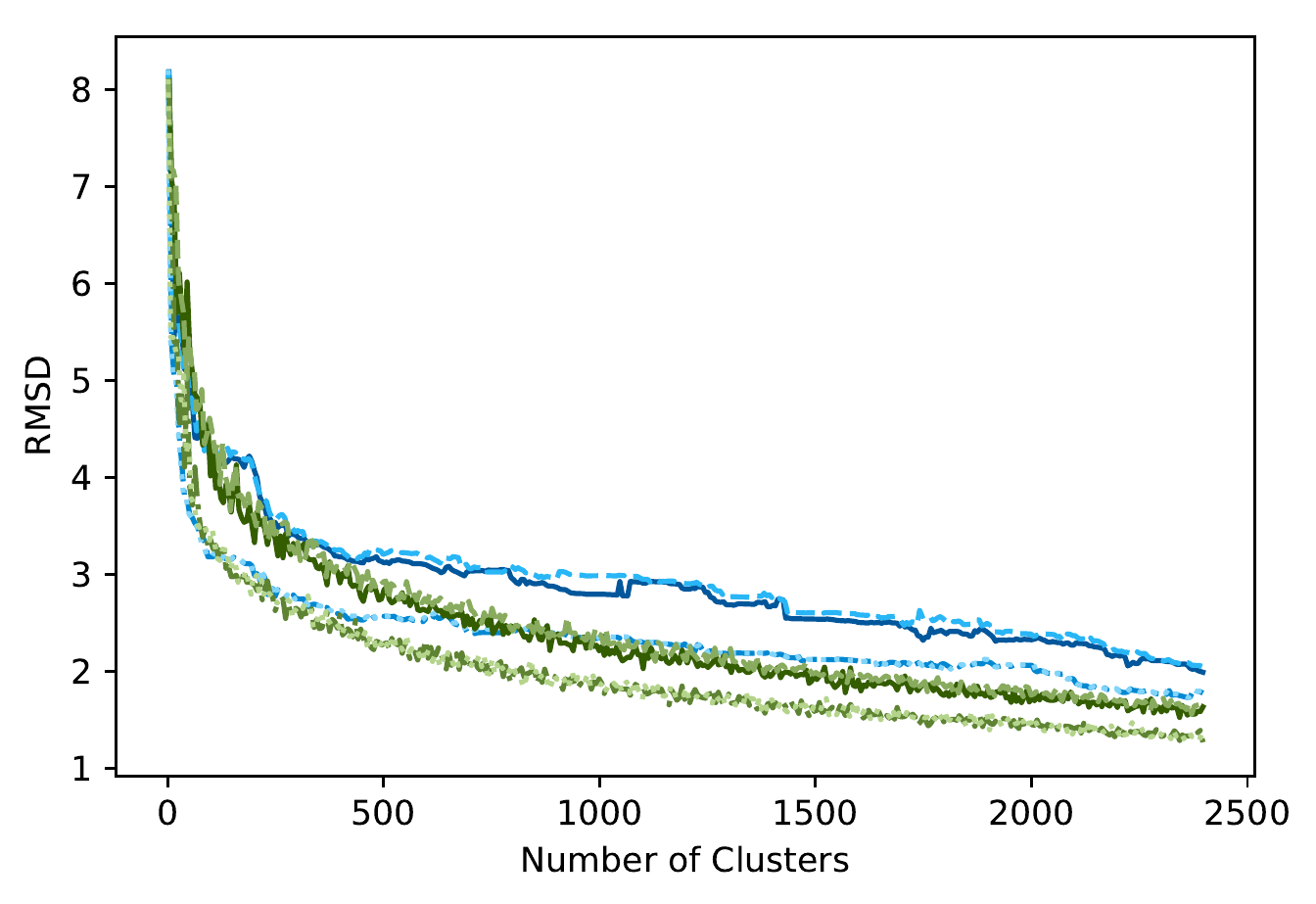}}
    \subfloat[Days]{\includegraphics[width=0.24\textwidth]{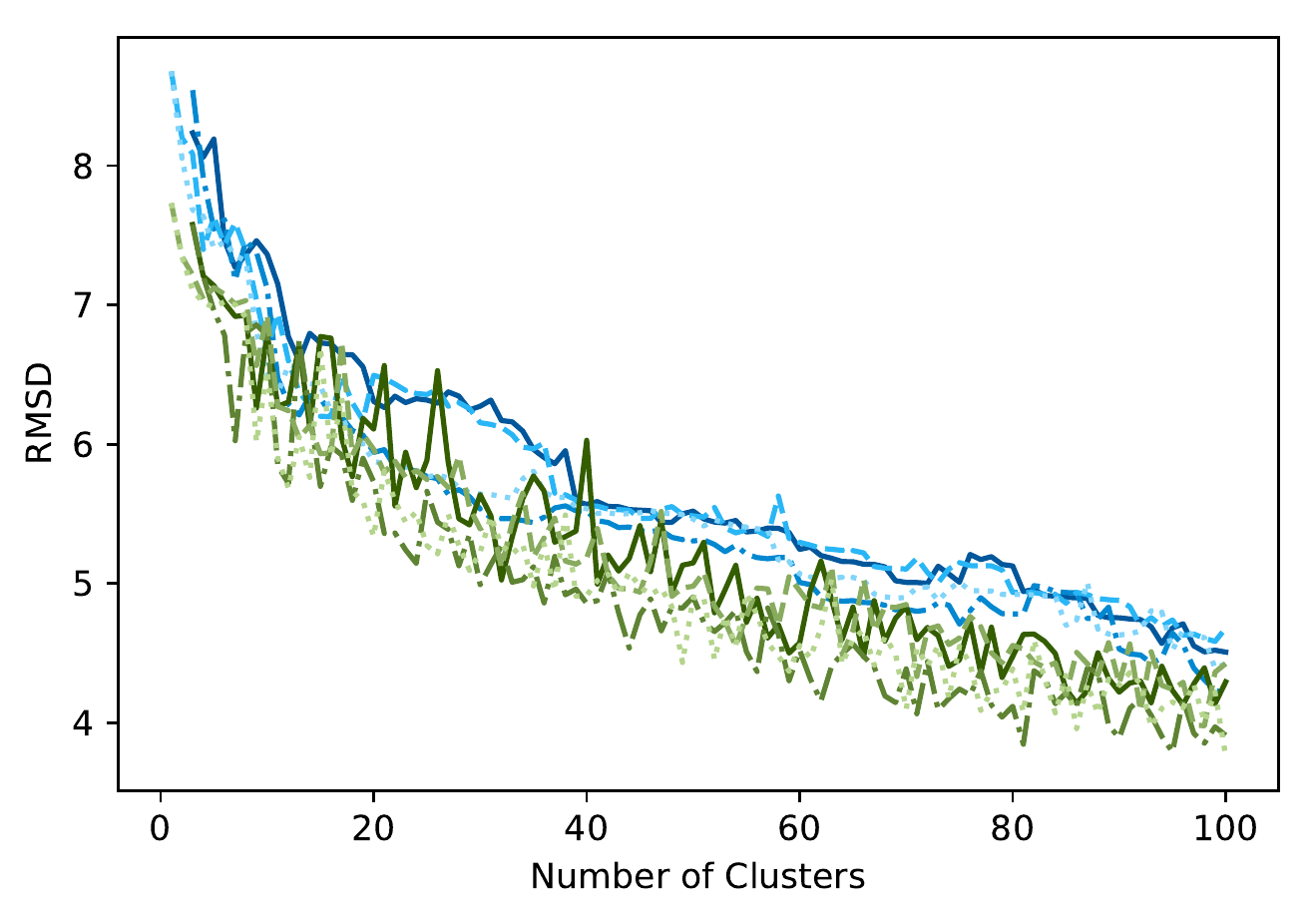}}
    \caption{RMSD of $CO_2$ Factor of Electricity}
    \label{co2}
\end{figure}

%\begin{figure}[h]
%    \centering
%    \includegraphics[width=0.48\textwidth]{RMSD_CO2_Factor_for_Electricity_days.pdf}
%    \caption{RMSD of $CO_2$ factor of electricity, Days}
%\end{figure}

For the $CO_2$ factor of electricity in Fig. \ref{co2}, the convergence rate is  much lower in the case of days than of hours. The decrease is almost linear, compared to exponential. In addition, there are high variations for days that are not present for hours. For 100 days, the RMSD is about 4.5 $gCO_2/kWh$ against 2 for an equivalent number of hours.

\begin{figure}[h]
    \centering
    \subfloat[Hours]{\includegraphics[width=0.24\textwidth]{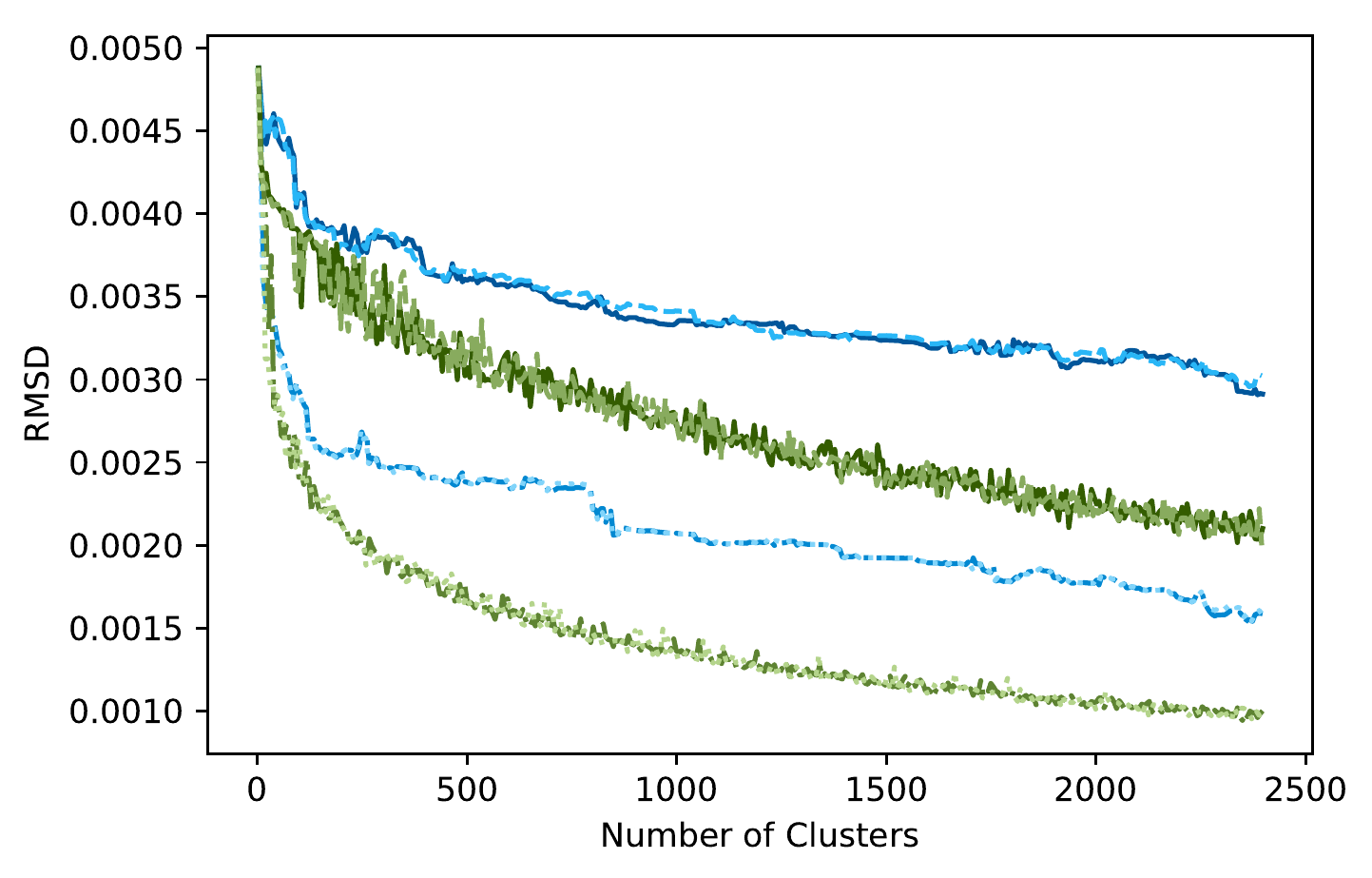}}
    \subfloat[Days]{\includegraphics[width=0.24\textwidth]{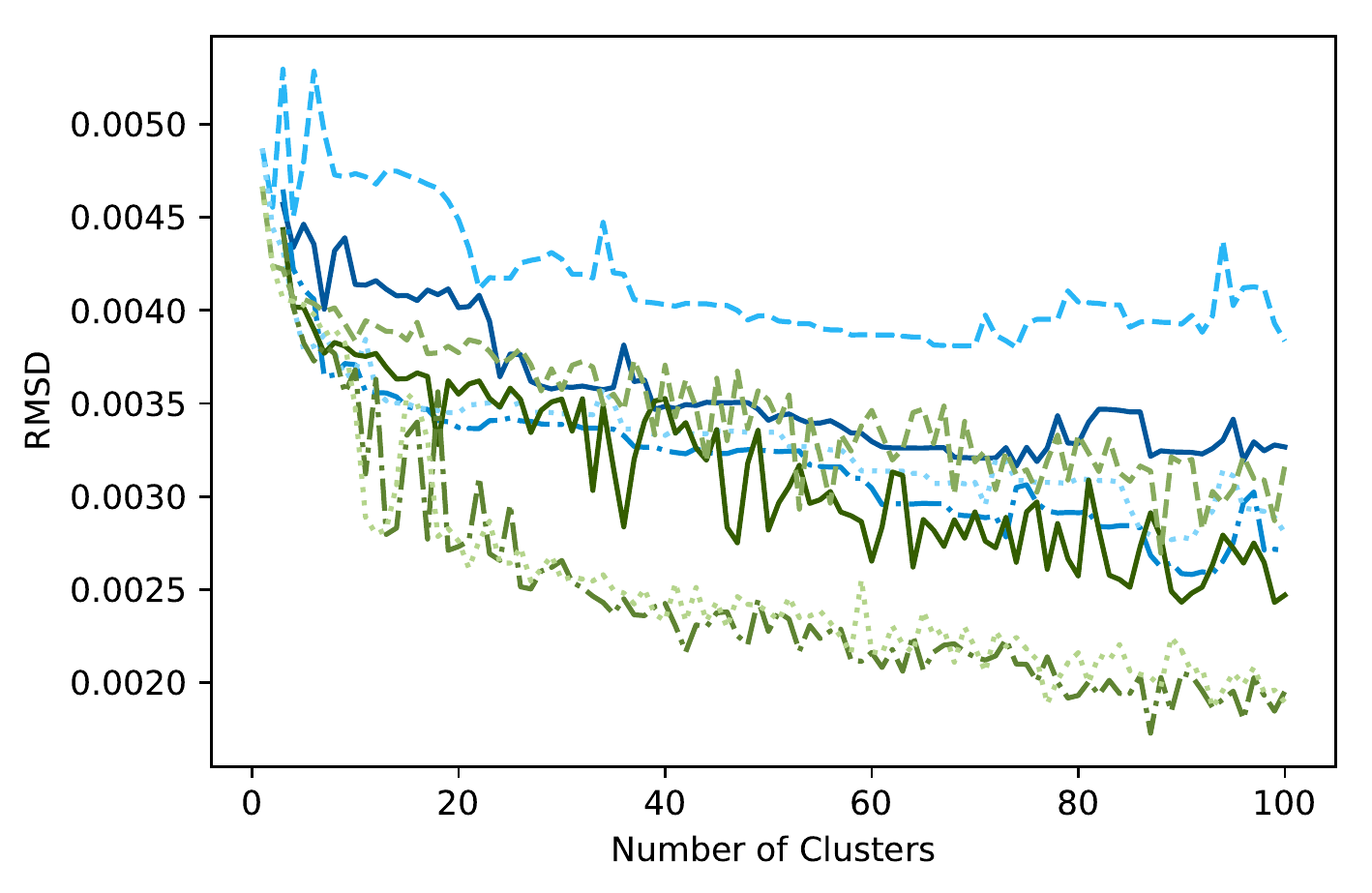}}
    \caption{RMSD of Spot Price}
    \label{spot}
\end{figure}

%\begin{figure}[h]
%    \centering
%    \includegraphics[width=0.48\textwidth]{RMSD_Spot_Price_days.pdf}
%    \caption{RMSD of Spot Price, Days}
%\end{figure}

In the case of spot price in Fig. \ref{spot}, the difference between the clustering performed using the standard deviation method and the range method seems very significant: for hours clustering between 0.001 and 0.0015$\euro/kWh$ or a factor of 2, the standard deviation is performing better. The difference between k-medoids and k-means is also considerably to the advantage of k-means: for hours clustering between 0.0005 and 0.001 $\euro/kWh$. For design days clustering the overall difference between methods is similar but there is more variability and some differences specific to this case. For instance, there are differences between the cases with and without heuristic, with the heuristic case performing better. Those differences are rather small for the standard deviation normalization and larger in the range case, especially in the k-medoids case.

%\begin{figure}[h]
%    \centering
%    \subfloat[Hours]{\includegraphics[width=0.24\textwidth]{RMSD_Temperature_hours.p%df}}
%    \subfloat[Days]{\includegraphics[width=0.24\textwidth]{RMSD_Temperature_days.pdf%}}
%    \caption{RMSD of Temperature}
%    \label{temperature}
%\end{figure}
%
%The errors for the temperature time series, Fig. \ref{temperature}, are very %similar to the overall ones commented before. The RMSD of temperature plateaus %rather quickly to around 2 for the hours, and 2.8 for the days.

The errors for the temperature time series are very similar to the overall ones commented before. The RMSD of temperature plateaus rather quickly to around 2 for the hours, and 2.8 for the days.

\begin{figure}[h]
    \centering
    \subfloat[Hours]{\includegraphics[width=0.24\textwidth]{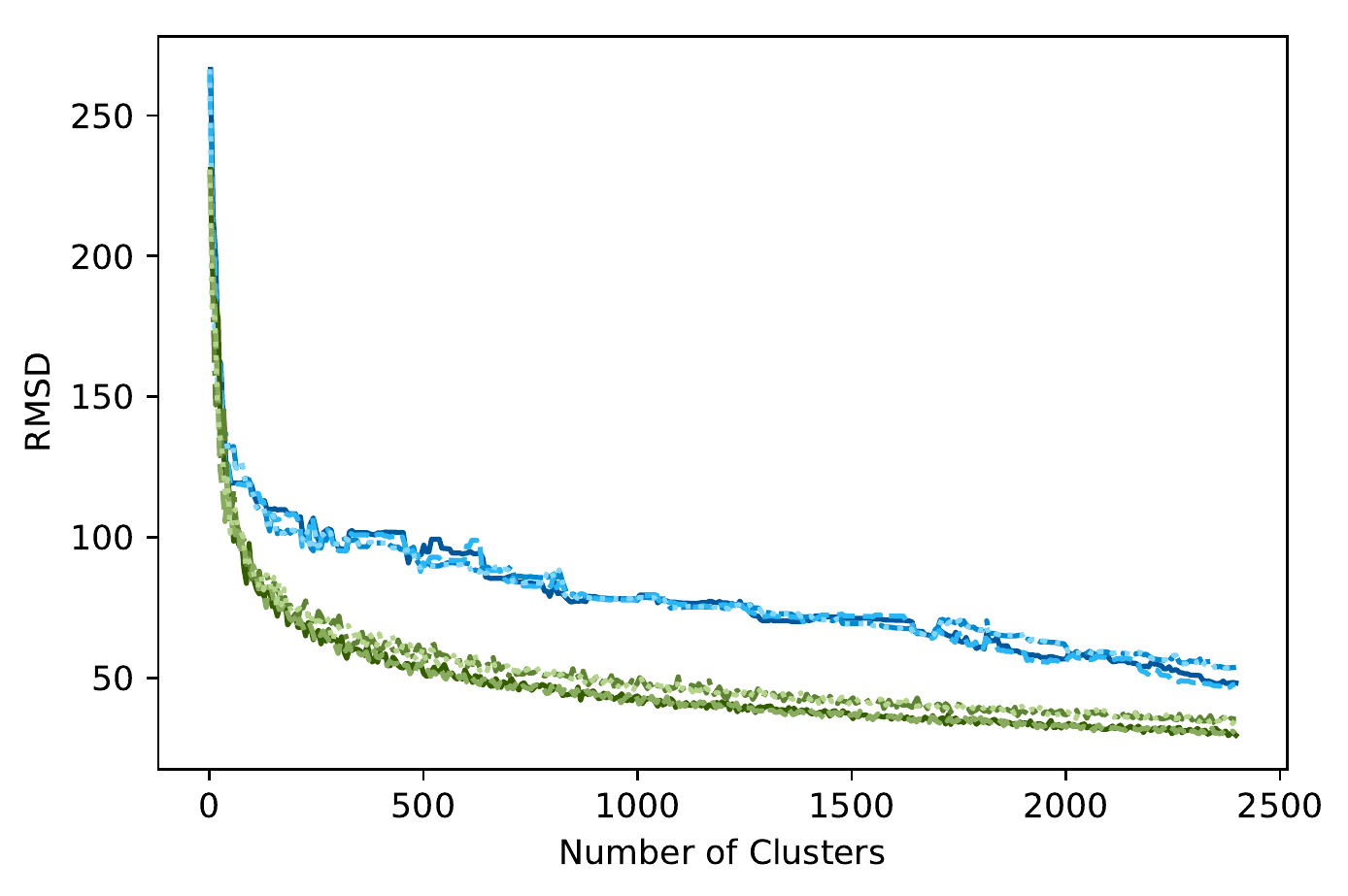}}
    \subfloat[Days]{\includegraphics[width=0.24\textwidth]{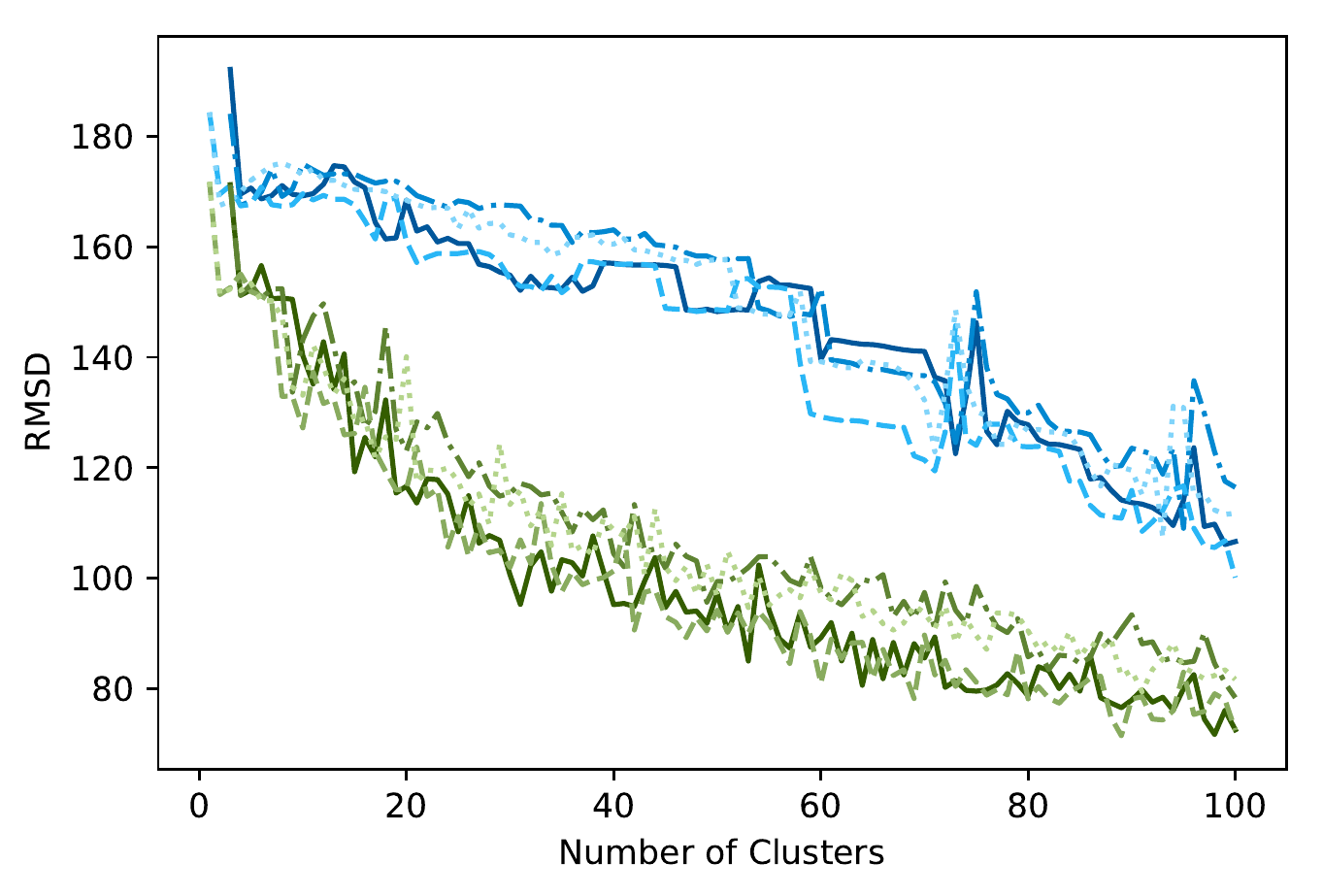}}
    \caption{RMSD of Irradiance on a Tilted Surface}
    \label{irr}
\end{figure}

%\begin{figure}[h]
%    \centering
%    \includegraphics[width=0.48\textwidth]{RMSD_Irr_tilt_days.pdf}
%    \caption{RMSD of Irradiance on a Tilted Surface, Days}
%\end{figure}

In the context of Zero Emission Neighborhoods, the irradiance has a very important role. Indeed, solar power is the main source of local (on the site) energy for neighborhoods. This means that solar irradiance and the production from the solar technologies will be crucial in compensating the emissions in the Zero Emission balance. Thus, in order to obtain designs that actually are Zero Emission, the precision of the clustering of the irradiance is essential.
The behaviour for the hours clustering, Fig. \ref{irr}, is similar to the overall behaviour. The RMSD for 100 clusters is around 35$W/m^2$ for k-means and 55$W/m^2$ for k-medoids. For days clustering, the convergence rate is slow and after 100 cluster, the RMSD is around 80 for k-means and 110 for k-medoids. The slow convergence rate means that for small numbers of clusters the difference between days and hours clustering is even worse. For 10 days, the RMSD is 140 for k-means and 170 for k-medoids. For 240 hours, the RMSD is about 60 for k-means and 100 for k-medoids. Those values are high in comparison to the standard test condition (STC) of solar panels of 1000 $W/m^2$.

\begin{figure}[h]
    \centering
    \subfloat[Hours]{\includegraphics[width=0.24\textwidth]{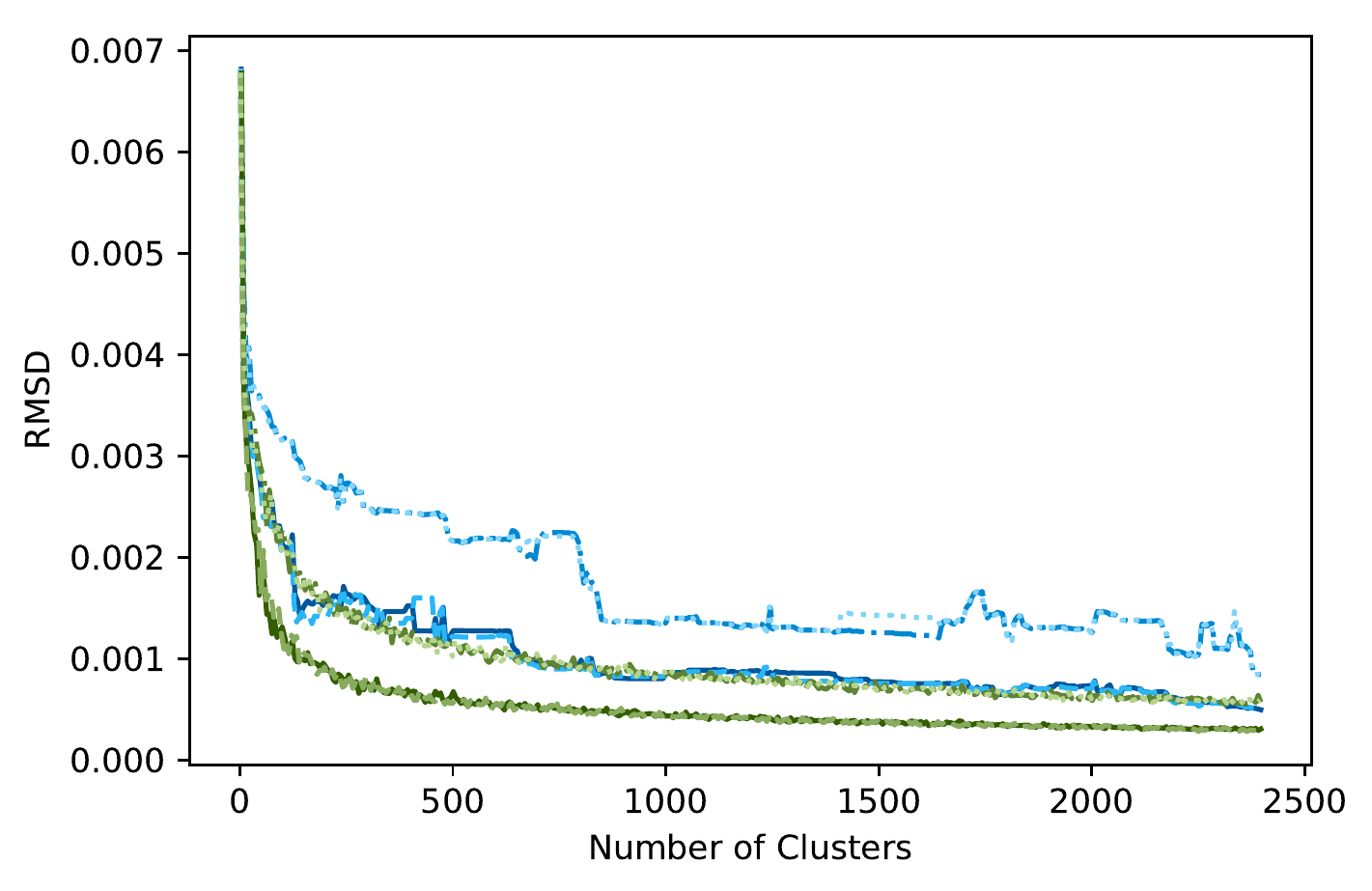}}
    \subfloat[Days]{\includegraphics[width=0.24\textwidth]{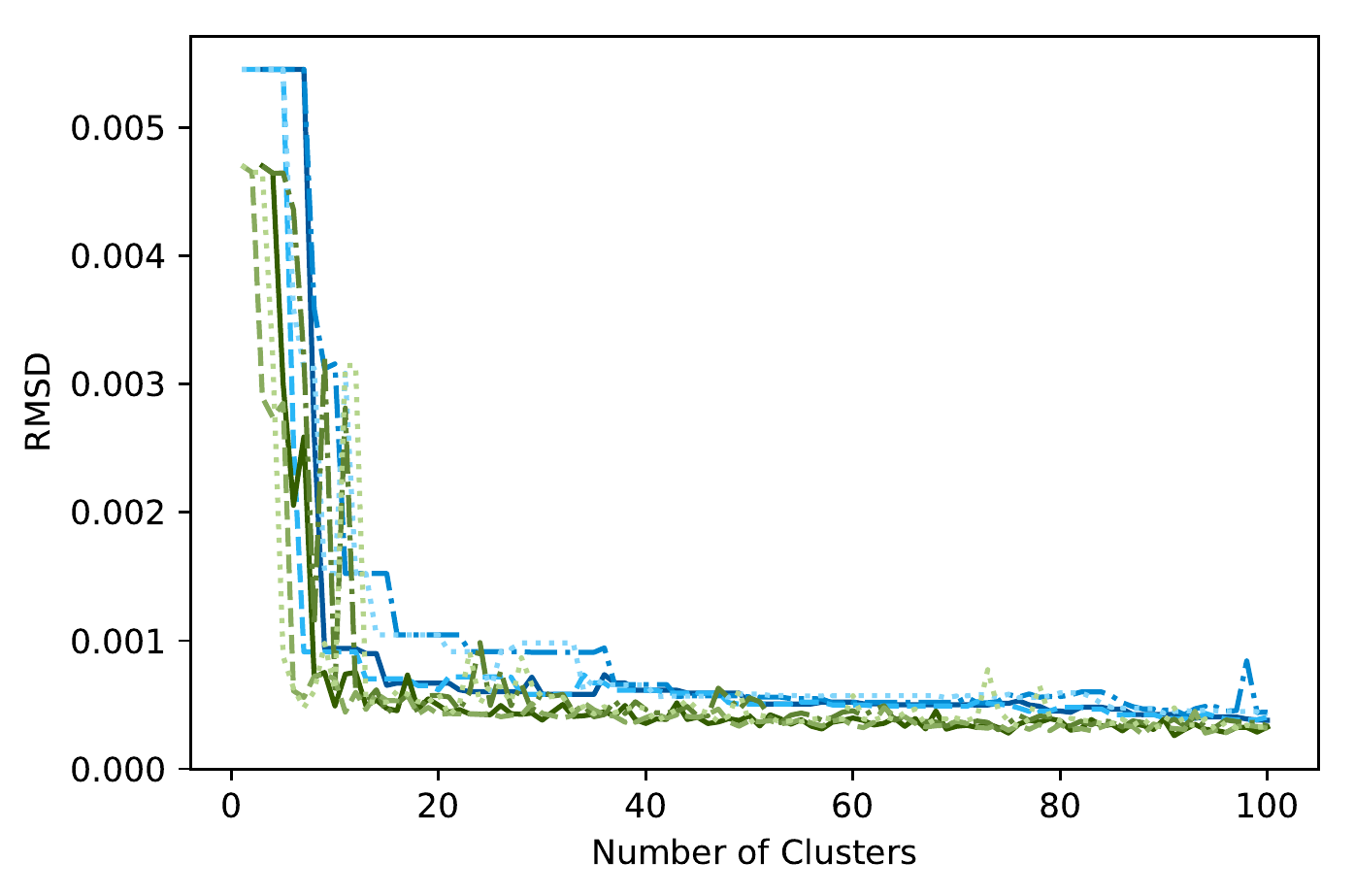}}
    \caption{RMSD of Electric Load in the Normal Offices}
    \label{elL}
\end{figure}

%\begin{figure}[h]
%    \centering
%    \includegraphics[width=0.48\textwidth]{RMSD_Electrical_Load_in_Normal_Offices_days.pdf}
%    \caption{RMSD of Electric Load in the Normal Offices, Days}
%\end{figure}

Only the performance for one of the three buildings is shown in this section. The other buildings can be found in the appendix.

For the electric load, Fig. \ref{elL}, in the case of days, the convergence has a steep rate but it happens slightly later around 10 days. After the convergence, the difference between all methods is close to zero. For the clustering on hours, the convergence is fast. The main difference from the behavior in the mean RMSD is that the cases with k-medoids and standard deviation normalization have a higher RMSD. The plateau is around 0.0013$Wh.m^{-2}.h^{-1}$ versus 0.0005$Wh.m^{-2}.h^{-1}$ for k-means range and 0.0008$Wh.m^{-2}.h^{-1}$ for the others.

%\begin{figure}[h]
%    \centering
%    \subfloat[Hours]{\includegraphics[width=0.24\textwidth]{RMSD_SH_Load_in_Normal_O%ffices_hours.pdf}}
%    \subfloat[Days]{\includegraphics[width=0.24\textwidth]{RMSD_SH_Load_in_Normal_Of%fices_days.pdf}}
%    \caption{RMSD of SH Load in the Normal Offices}
%\end{figure}

%\begin{figure}[h]
%    \centering
%    \includegraphics[width=0.48\textwidth]{RMSD_SH_Load_in_Normal_Offices_days.pdf}
%    \caption{RMSD of SH Load in the Normal Offices, Days}
%\end{figure}

%\begin{figure}[h]
%    \centering
%    \subfloat[Hours]{\includegraphics[width=0.24\textwidth]{RMSD_DHW_Load_in_Normal_%Offices_hours.pdf}}
%    \subfloat[Days]{\includegraphics[width=0.24\textwidth]{RMSD_DHW_Load_in_Normal_O%ffices_days.pdf}}
%    \caption{RMSD of DHW Load in the Normal Offices}
%\end{figure}

The RMSD for the SH and DHW time series behave as the mean of the RMSDs. The mean of the RMSD is influenced greatly by the loads because they behave similarly and because of the presence of 3 time series for each building.

%\begin{figure}[h]
%    \centering
%    \includegraphics[width=0.48\textwidth]{RMSD_DHW_Load_in_Normal_Offices_days.pdf}
%    \caption{RMSD of DHW Load in the Normal Offices, Days}
%\end{figure}

%individual results: irr_tilt, co2, ...

Another metric of interest is the Yearly Average Error (YAE), this metric allows us to have information about the distribution of the error. With RMSD, there is no information on the sign of the errors. YAE allows to know if the errors are, on average, compensated or rather cumulate from timestep to timestep.

\begin{table}[t]
\renewcommand{\arraystretch}{1.15}
    \centering
    \caption{Yearly Average Error (YAE) and RMSD for 10 and 100 days and equivalent number of hours for the irradiance with k-medoids (STD.:Standard Devation, H: with heuristic, \cancel{H}: without heuristic)}
    \begin{tabular}{@{} l r r r r r r r r @{}}
    \hline
         & \multicolumn{4}{c}{\textbf{Days}} & \multicolumn{4}{c}{\textbf{Hours}}\\
            \hline
        & \multicolumn{2}{c}{\textbf{STD.}} & \multicolumn{2}{c}{\textbf{Range}} & \multicolumn{2}{c}{\textbf{STD.}} & \multicolumn{2}{c}{\textbf{Range}}\\
            \hline
        & \textbf{H} & \textbf{\cancel{H}} & \textbf{H} & \textbf{\cancel{H}} & \textbf{H} & \textbf{\cancel{H}} & \textbf{H} & \textbf{\cancel{H}}\\
        \hline
        YAE 10 & -34 & -24 & -37 & -38 & -18 & -18 & -13 & -15 \\
        YAE 100 & -3.5 & -0.67 & -0.26 & -2.3 & -3.1 & -2.8 & -3.0 & -2.1 \\
        RMSD 10 & 175 & 173 & 169 & 170 & 95.3 & 95.6 & 107 & 107 \\
        RMSD 100 & 116 & 112 & 107 & 100 & 53.9 & 53.8 & 48.3 & 47.8 \\
         \hline
    \end{tabular}
    \label{tab:cluster_res}
\end{table}

The results for k-means are not in Table \ref{tab:cluster_res} because the YAE stays at 0 for all number of clusters and all cases. Instead the values of the RMSD are presented in Table \ref{tab:cluster_res2}.

\begin{table}[t]
\renewcommand{\arraystretch}{1.15}
    \centering
    \caption{RMSD for 10 and 100 days and equivalent number of hours for the irradiance with k-means (STD.:Standard Deviation, H: with heuristic, \cancel{H}: without heuristic)}
    \begin{tabular}{@{} l r r r r r r r r @{}}
    \hline
         & \multicolumn{4}{c}{\textbf{Days}} & \multicolumn{4}{c}{\textbf{Hours}}\\
            \hline
        & \multicolumn{2}{c}{\textbf{STD.}} & \multicolumn{2}{c}{\textbf{Range}} & \multicolumn{2}{c}{\textbf{STD.}} & \multicolumn{2}{c}{\textbf{Range}}\\
            \hline
        & \textbf{H} & \textbf{\cancel{H}} & \textbf{H} & \textbf{\cancel{H}} & \textbf{H} & \textbf{\cancel{H}} & \textbf{H} & \textbf{\cancel{H}}\\
        \hline
        RMSD 10 & 143 & 133 & 140 & 127 & 71.0 & 72.2 & 63.7 & 66.6 \\
        RMSD 100 & 78.4 & 81.7 & 72.6 & 72.5 & 35.5 & 35.1 & 29.7 & 30.5 \\
         \hline
    \end{tabular}
    \label{tab:cluster_res2}
\end{table}

Comparing the RMSD and YAE from Table \ref{tab:cluster_res} and Table \ref{tab:cluster_res2} gives us insights in how much the errors in irradiance cancel each other, at least in terms of annual values. 
In the case of irradiance, the negative signs first informs us that it is under-represented. The difference between the RMSD and the YAE values also suggest that the errors tends to be compensated by one another and they compensate completely in the k-means cases. 
In general, the hours clustering performs better than the daily one. k-means is better than k-medoids in terms of YAE for the same reasons that it is better for RMSD. 
The performance of STD or range on their own or in addition to heuristic is not consistent but the gains here are less big than between days and hours clustering.

From the results presented in this section, k-means and hours clustering are the best choices. For instance, with a focus on the irradiance, the choice would be range and heuristic. Overall the biggest impact can be made by choosing the correct clustering algorithm and the correct resolution. When it comes to the normalization method and the use of heuristic, the choice has less importance and varies depending on the goal. However there appears to be better results with the range normalization and without the heuristic.
These results are however not enough. They only display some metrics for how close the clusters come to the original data. This does not guarantee that the one performing best in this section would also perform best in the optimization.

\section{Models and Implementation}

In this section, the main equations of the ZENIT model are presented along with two variations for using either representative days or hours then the implementation and data used is briefly presented. The variations will be called M0 and M1 and are based on \cite{GABRIELLI2018408}.

ZENIT aims is to design the energy system of a neighborhood so that it can be Zero Emission during its lifetime. Thus, it considers the investment as well as the operation of the neighborhood to find the cost optimal solution.
The objective function is:
\textit{Minimize}:
\begin{multline}
    b^{HG}\cdot C^{HG} + \sum_{b}\sum_{i}  \Big( (C_{i,b}^{var,disc} + \frac{C_{i,b}^{maint}}{\varepsilon^{tot}_{r,D}})\cdot x_{i,b} +\\ C_{i,b}^{fix,disc} \cdot b_{i,b} \Big) 
     + \sum_{t_\kappa}\frac{\sigma(\kappa)}{\varepsilon^{tot}_{r,D}}\Big(\sum_b\sum_{f}f_{f,t,b} \cdot P_{f}^{fuel}\\
     + (P_{t}^{spot}+P^{grid}+P^{ret})\cdot (y_{t}^{imp} \\
     +\sum_b\sum_{est}y_{t,est,b}^{grid\_imp,bat})-P_{t}^{spot}\cdot y_{t}^{exp}\Big)
\end{multline}

It considers the investment cost in technologies ($C_{i,b}^{var,disc}$, $C_{i,b}^{fix,disc}$) and the heating grid ($C^{HG}$), as well as operation and maintenance related costs ($C_{i,b}^{maint}$). A binary variable controls the investment in the heating grid ($b^{HG}$). The subscript used in the equations are $b$ for the buildings, $i$ for the technologies, $t$ for the timesteps, $f$ for fuels and $est$ for batteries. $\varepsilon$ are the discount factors with interest rate $r$ for the duration of the study $D$. $x_{i,b}$ is the capacity of the technologies and $b_{i,b}$ the binary related to whether it is invested in or not. $\sigma(\kappa)$ is the number of occurrences of cluster $\kappa$ in the full year and $t_\kappa$ is the timestep in the cluster. $P$ are the prices of fuel, electricity on the spot market, grid tariff or retailer tariff. $f$ is the consumption of fuel and $y$ are the imports or exports of electricity.

In order to fulfill the Zero Emission requirement, the following constraint, called the Zero Emission Balance is necessary:

\begin{multline}
    \phi^{CO_2}_{e,t}\sum_{t_\kappa}\sigma(\kappa)\Big(y_{t}^{imp}+\sum_b\sum_{est}y_{t,est,b}^{grid\_imp,bat}\Big) \\
    + \sum_{t_\kappa}\sigma(\kappa)\sum_b\sum_{f} \phi^{CO_2}_f \cdot f_{f,t,b} \leq  \phi^{CO_2}_{e,t}\cdot\sum_{t_\kappa}\sigma(\kappa) \\
    \bigg( \sum_b\sum_{est}\eta_{est}\cdot\Big(\alpha_{ZEN}\cdot y_{t,est,b}^{grid\_exp,bat}+y_{t,est,b}^{prod\_exp,bat}\Big) \\
    +\sum_b\sum_{g}\Big(y_{t,g,b}^{exp}\Big)\bigg)
\end{multline}

It forces the emissions of $CO_2$ to be at least equal to the compensations. The principle of the compensation is that the energy produced in the neighborhood, by renewable sources, that is exported to the national grid reduces the global production. The corresponding amount of saved $CO_2$ is counted as compensation for the neighborhood. The $CO_2$ factors are represented by $\phi^{CO_2}_{e,t}$ for electricity and $\phi^{CO_2}_f$ for other fuels. $\eta_{est}$ is the charging efficiency of the battery.

Other equations include load balances for electricity, DHW and SH. They require the production and import to be equal to the consumption and exports.

In ZENIT, different scales of technologies can be invested in. There are technologies at the building level and technologies at the neighborhood level gathered at a "Production Plant".

The investment in technologies are bounded, from above and below or only from above. When it is bounded from below, the capacity of that technology is defined as a semi-continuous variable thanks to a binary. In addition, for technologies at the plant level, the heating grid needs to be there:
$\forall i$ 
\begin{equation}
     x_{i,'Production Plant'} \leq X_{i}^{max} \cdot b^{HG}
\end{equation}
$X_{i}^{max}$ is the maximum capacity of technology $i$.
For the heat pumps in the buildings, the production and consumption are defined as follow:
\begin{equation}\label{eq:copsh}
    d_{hp,b,t}^{SH}=\frac{q_{hp,b,t}^{SH}}{COP_{hp,b,t}^{SH}}
\end{equation}
\begin{equation}\label{eq:copdhw}
    d_{hp,b,t}^{DHW}=\frac{q_{hp,b,t}^{DHW}}{COP_{hp,b,t}^{DHW}}
\end{equation}
\begin{equation}\label{eq:HP_lim}
    \frac{d_{hp,b,t}^{DHW}}{P^{input,max,DHW}_{hp,b,t}} + \frac{d_{hp,b,t}^{SH}}{P^{input,max,SH}_{hp,b,t}} \leq x_{hp,b}
\end{equation}

Equations \ref{eq:copsh} and \ref{eq:copdhw} link the heat produced to the COP and the electrical consumption. The COPs are different for SH and DHW due to different temperature set points. They also depend on the outside temperature and are calculated before the optimization.
Equation \ref{eq:HP_lim} regulates how the heat pump can be used for both SH and DHW and enforces that the capacity invested is not exceeded. $P^{input,max}$ represents the maximum power input to the heat pump at the timestep based on the temperature set point for a 1kW unit. $d_{hp,b,t}^{SH}$ and $d_{hp,b,t}^{SH}$ represent the electric consumption of the heat pump for SH and DHW while $q_{hp,b,t}^{DHW}$ and $q_{hp,b,t}^{DHW}$ are the heat production.

Another binary variable is used for part load limitations. This binary concerns the operation and is defined for every timestep for each relevant technology, which can lead to a large number of binary variable. No minimum up- or downtime is used. 

The handling of the storages is what differentiates model M0 and M1. Model M0 is not able to handle seasonal storage while model M1 can be used for that.
In ZENIT, each battery is modelled as 2 separate virtual batteries, with one connecting the neighborhood to the grid: allowing import and export between the grid and the battery, and import from the battery to the neighborhood's loads, and the other connecting the technologies producing electricity in the neighborhood and the neighborhood's loads: allowing exports to the grid and to the neighborhood. This distinction provides traceability of the electricity in the batteries. The origin of the electricity is important because of the different $CO_2$ factors.

Both models are presented in the following subsections.

\subsection{Model M0}

Model M0 uses a classical formulation for storages in models using clustering and that do not need seasonal storage. The equation for electric and heat storages are similar, so only a generic equation is presented, fitting both cases.

$\forall \kappa, t\in [1,T^{clu}],st,b$
\begin{equation}
    v^{stor}_{\kappa,t,st}=v^{stor}_{t-1,st}+\eta^{stor}_{st}\cdot q_{t,st}^{ch} -q_{t,st}^{dch}
\end{equation}

$\forall t\in [0,T^{clu}],st,b$ 
\begin{equation}
     v^{stor}_{\kappa,t,st,b} \leq x_{st,b}
\end{equation}

\begin{multicols}{2}\noindent
    \begin{equation}
     q^{ch}_{\kappa,t,st,b} \leq \dot{Q}_{st}^{max}
    \end{equation}\noindent
    \begin{equation}
     q^{dch}_{\kappa,t,st,b} \leq \dot{Q}_{st}^{max}
    \end{equation}
\end{multicols}
$\forall p,st,b,\kappa$
\begin{equation}\label{eq:endval}
    v^{stor}_{\kappa,0,st,b}=v^{stor}_{\kappa,T^{clu},st,b}
\end{equation}

The state of charge of the storage $st$ (either heat or electric storage) is represented by $v^{stor}$ while $q{ch}$ and $q{dch}$ are the energy charged and discharged. The maximum charge and discharge rate is ${Q}_{st}^{max}$.
The differences between this model and a model with full year data is that the starting value of the storage is "free" at the beginning of each cluster instead of only at the beginning of the year. This model is valid for different $T^{clu}$ even though in our case it is 24 for days clustering. It allows for a daily operation of the storage. Different values of $T^{clu}$ could be used for allowing different ranges of operation of the storage. A bigger value allows longer operation but probably increases the number of timesteps to get the same clustering precision and reducing it reduces the possible range. The daily range makes sense because of the daily cycle of the loads, that allows us to make the assumption of equation \ref{eq:endval}.
This model does not make sense with $T^{clu}=1$, i.e. the hours clustering, because the resolution of the data used is also one hour; hourly storage operation does not make sense.

\subsection{Model M1}

In model M1, the main difference with model M0 is that the storage level equation becomes:
$\forall \kappa, t\in [1,8760],st,b$
\begin{equation}
    v^{stor}_{\kappa,t,st}=v^{stor}_{t-1,st}+\eta^{stor}_{st}\cdot q_{t_{\kappa},st}^{ch} -q_{t_{\kappa},st}^{dch}
\end{equation}
The end value of the storage constraint is also replaced by:
$\forall st,b$
\begin{equation}
    v^{stor}_{0,st,b}=v^{stor}_{8760,st,b}
\end{equation}

Where $t_{\kappa}$ is the time corresponding to t in the cluster. It is found by using the sequence of cluster's representatives ($\xi$) either directly for the hourly case or with the day number corresponding to t and the hour in the day :

\noindent
\begin{minipage}[c]{0.29\linewidth} 
    \textbf{Hours}:
\end{minipage} % no space if you would like to put them side by side
\begin{minipage}[c]{0.7\linewidth}
    \textbf{Days}:
\end{minipage}
\begin{minipage}[c]{0.29\linewidth}
    \begin{equation} \label{hour}
        t_{\kappa} = \xi(t)
    \end{equation}
\end{minipage} % no space if you would like to put them side by side
\begin{minipage}[c]{0.7\linewidth}
    \begin{equation} \label{day}
        t_{\kappa} = \xi\bigg(\left\lfloor\frac{t}{24}\right\rfloor \bigg) + t-\left\lfloor\frac{t}{24}\right\rfloor \cdot 24
    \end{equation}
\end{minipage}

This means that the storage level is not decoupled between the different clusters. The charging and discharging is defined for each timestep in each cluster but the storage level is defined for every hour in the year. 
This model comes from the assumption that days or hours with similar conditions in terms of the time series (loads, spot price, temperature,...), i.e. belonging to one cluster will behave in the same way in terms of charging and discharging of the storage.
This formulation however comes at the expense of longer computation time.
Both hourly and days clustering can be used with M1.

In \cite{GABRIELLI2018408}, another variation is presented to improve further model M1 by defining only the variables related to operations binary variables (on/off status) for each cluster while other variables are defined for each hour of the year. That model has not been implemented because it increases the computation time even more.

\subsection{Implementation}

The model is implemented on a test case based on a small neighborhood, a campus at Evenstad in Norway. The buildings are gathered in three categories to only have three buildings in the optimization. We assume every building has a hydronic system.

The economical and technical data of the technologies are taken from the Danish Energy Agency\footnote{\url{https://ens.dk/en/our-services/projections-and-models/technology-data}}. In total, 22 technologies are implemented with, at the building level: solar panel, solar thermal, air-air heat pump, air-water heat-pump, ground source heat pump, bio boiler with wood logs or pellets, electric heater and electric boiler, biomethane boiler, biogas and biomethane CHP; and at the neighborhood level: biogas boiler, wood chips and pellets boiler and CHPs, ground source heat pump and electric boiler.

The spot price of electricity is obtained from Nordpool's website \footnote{\url{https://www.nordpoolgroup.com/Market-data1/#/nordic/table}}.
The temperature data comes from Agrometeorology Norway\footnote{\url{https://lmt.nibio.no}, F\aa vang station}. The solar irradiance (diffuse horizontal (DHI) and direct normal(DNI)) are obtained from Solcast \footnote{\url{https://solcast.com.au}}. The irradiance on a tilted surface $IRR^{Tilt}$ which is an input of the clustering is derived from the DHI and DNI with:

\begin{multline}\label{eq:irr}
    IRR^{Tilt}_{t}=DHI_{t}\frac{1+cos(\phi_1)}{2} \\
    + \alpha \cdot \Big(DNI_{t}+DHI_{t}\Big)\frac{1-cos(\phi_1)}{2} \\
    + DNI_{t}\bigg( \frac{cos(\varphi_{t}) \cdot sin(\phi_1) \cdot cos(\phi_2-\psi_{t})}{sin(\varphi_{t})}\\+\frac{sin(\varphi_{t}) \cdot cos(\phi_1)}{sin(\varphi_{t})}\bigg)
\end{multline}

We assume that for some sun positions (sun elevations ($\varphi$) below 1 degree and sun azimuths ($\psi$) between -90 and 90 degrees), no direct beam reaches the panels. This means that the last term of equation \ref{eq:irr} is removed at such times.
We use a constant albedo factor ($\alpha$) of 0.3 for the whole year. Hourly albedo values could also be used to reflect the impact of snow in the winter better.
The tilt angle of the solar panel is $\phi_1$; the orientation of the solar panel regarding the azimuth is $\phi_2$.

The hourly $CO_2$ factors of electricity are obtained with the methodology presented in \cite{clauss19} while the other $CO_2$ factors come from \cite{clark13}.

The prices of wood pellets comes from \cite{nobio}, the price of wood logs from \cite{eubionet3}, the price of wood chips from \cite{tromborg15} and the price of biogas from \cite{eba16}.

The electric and heat load profiles for the campus are derived from \cite{lindberg_impact_2017}.
The domestic hot water (DHW) and Space Heating (SH) are then based on the time series from a passive building in Finland \cite{pal_energy_2016}.

The model is implemented in Python and is solved using Gurobi. It is run on a laptop  with an Intel Core i7-7600U dual core processor at 2.8Ghz and 16GB of RAM.

\section{Model Results and Discussion}

In this section we present the results obtained with the different clustering methods and variations from the earlier sections. We always use the heuristic in order to guarantee that the peak load is covered.

\subsection{Simplified Model}
In order to get a reference objective value to base our analysis on, a simplified version of the model is run. This simplified model leaves out several of the constraints using binaries, namely the part load constraints, the minimum investment capacity (turning the semi-continuous variables into continuous variables) and changing the cost function from $a\cdot x+b$ to $c\cdot x$. Without simplifying the model, solving the model with 365 days or 8760 hours would take too long (several days). It is important to note that this simplified model is not directly obtained by removing constraints but by setting the input associated to the binary to zero. For example, the fixed investment costs and the minimum capacity are set to 0 but the constraints are still there. In the case of the minimum load during operation, the minimum loads are set to 0 but the related constraints are not written when the model is generated in Gurobi.
The results for the non-simplified model are presented after without a reference value.

Because M1 allows for seasonal storage modelling while M0 does not and in order to obtain results that can be compared more easily between M0 and M1, the storages at the neighborhood level were not included in the technological option input in this study.

We chose the number of days and hours in this section graphically at the elbow of the curves in Fig. \ref{NRMSD}. The number of clusters is chosen so that adding clusters does not bring considerable improvements. For the case of hours, this corresponds to around 120 hours; 96 and 144 hours are also studied as a 20\% variation. We also consider the corresponding number of days, i.e. 4,5,6. Indeed this gives an equal number of timesteps in the optimization but the performance of the clustering on days for such low numbers of days should give poor result considering Fig. \ref{NRMSD}. 
In addition we choose a number of design days with similar graphical elbow considerations. However, we consider Fig. \ref{irr} instead of the NRMSD figure because in the case of clustering days the performances for the irradiance were converging slower. This leads us to choose 30 days. We also take the 20\% variations, which corresponds to 24 and 36 days.

\begin{table}[t]
\renewcommand{\arraystretch}{1.15}
    \centering
    \caption{Variations in objective value from the reference for different numbers of representative days for M0 with simplified model (\textbf{STD}: Standard Devation, \textbf{R}:Range), Reference Value for 365 days: \textbf{2,056,849} $\euro$}
    \begin{tabular}{@{} l r r r r r r r  @{}}
    \hline
          & & \multicolumn{6}{c}{\textbf{Days}}\\
            \hline
       % & & \multicolumn{2}{c}{4} & \multicolumn{2}{c}{5} & \multicolumn{2}{c}{6} & \multicolumn{2}{c}{24} & \multicolumn{2}{c}{30} & \multicolumn{2}{c}{36} & \multicolumn{2}{c}{96} & \multicolumn{2}{c}{120} & \multicolumn{2}{c}{144}\\
         & &  \textbf{4} & \textbf{5} & \textbf{6} & \textbf{24} & \textbf{30} & \textbf{36}\\
     %   \hline
    %    &  & STD & H & STD & H & STD & H & STD & H & STD & H & STD & H & STD & H\\
        \hline
        \multirow{2}{*}{\textbf{k-means}}  & STD & -10.29 & -9.50 & -9.42 & -6.14 & -5.21 & -4.82   \\
                                          & R & -10.29 & -10.64 & -9.68 & -4.80 & -4.74 & -3.82  \\
        \multirow{2}{*}{\textbf{k-medoids}} & STD & 28.27 & 22.71 & 33.53 & 9.61 & 10.04 & 7.49 \\
                                           & R & 11.57 & 8.78 & 23.16 & 5.36 & 4.84 & 8.27 \\
        \hline
    \end{tabular}
    \label{res_simplM0}
\end{table}

From Table \ref{res_simplM0}, k-means range seems to be the overall best choice, but it underestimates the objective value. k-medoids constantly overestimates the objective value, with significant errors for low numbers of days. On the other hand, k-means gives good results even for a low number of days. %Assuming the same behaviour as in the previous section, we can also note a better performance of the range method.

%\begin{table}[t]
%\renewcommand{\arraystretch}{1.15}
%    \centering
%    \caption{Variations in objective value from the reference for different number of representative days and %hours for M1(\textbf{STD}: Whiten, \textbf{R}:Range), Reference Value for a complete year: \textbf{2,060,612} %$\euro$ }
%    \begin{tabular}{@{} l r r r r r r  r @{}}
%    \hline
%          & & \multicolumn{3}{c}{\textbf{Days}} & \multicolumn{3}{c}{\textbf{Hours}}\\
%            \hline
%       % & & \multicolumn{2}{c}{4} & \multicolumn{2}{c}{5} & \multicolumn{2}{c}{6} & \multicolumn{2}{c}{24} & %\multicolumn{2}{c}{30} & \multicolumn{2}{c}{36} & \multicolumn{2}{c}{96} & \multicolumn{2}{c}{120} & %\multicolumn{2}{c}{144}\\
%        & &  \textbf{4} & \textbf{5} & \textbf{6} & \textbf{96} & \textbf{120} & \textbf{144}\\
%     %   \hline
%    %    &  & STD & H & STD & H & STD & H & STD & H & STD & H & STD & H & STD & H\\
%        \hline
%        \multirow{2}{*}{\textbf{k-means}} & STD & -10.16 & -9.18 & -8.78 & -5.58 & -5.54 & -4.95 \\
%                                         & R &  &  &  &  &  &  \\
%         \multirow{2}{*}{\textbf{k-medoids}} & STD & 28.42 & 22.80 & 33.55 & 8.45 & 11.06 & 9.73 \\
%                                            & R &  &  &  &  &  &  \\
%         \hline
%    \end{tabular}
%    \label{res_simpl}
%\end{table}

\begin{table}[t]
\renewcommand{\arraystretch}{1.15}
    \centering
    \caption{Variations in objective value from the reference for different number of representative days for M1 with simplified model (\textbf{STD}: Standard Devation, \textbf{R}:Range), Reference Value for a complete year: \textbf{2,060,612} $\euro$ }
    \begin{tabular}{@{} l r r r r r r  r @{}}
    \hline
       % & & \multicolumn{2}{c}{4} & \multicolumn{2}{c}{5} & \multicolumn{2}{c}{6} & \multicolumn{2}{c}{24} & \multicolumn{2}{c}{30} & \multicolumn{2}{c}{36} & \multicolumn{2}{c}{96} & \multicolumn{2}{c}{120} & \multicolumn{2}{c}{144}\\
        & &  \textbf{4} & \textbf{5} & \textbf{6} & \textbf{24} & \textbf{30} & \textbf{36}\\
     %   \hline
    %    &  & STD & H & STD & H & STD & H & STD & H & STD & H & STD & H & STD & H\\
        \hline
        \multirow{2}{*}{\textbf{k-means}} & STD & -10.16 & -9.18 & -8.78 & -6.07 & -5.38 & -6.14 \\
                                         & R & -10.16 & -10.38 & -9.09 & -5.03 & -5.38 & -4.44
                                         \\
         \multirow{2}{*}{\textbf{k-medoids}} & STD & 28.42 & 22.80 & 33.55 & 9.64 & 10.08 & 7.54 \\
                                            & R & 11.78 & 8.91 & 23.19 & 5.40 & 4.90 & 8.31 \\
         \hline
    \end{tabular}
    \label{res_simplM1_1}
\end{table}

\begin{table}[t]
\renewcommand{\arraystretch}{1.15}
    \centering
    \caption{Variations in objective value from the reference for different number of representative hours for M1 with simplified model (\textbf{STD}: Standard Devation, \textbf{R}:Range), Reference Value for a complete year: \textbf{2,060,612} $\euro$ }
    \begin{tabular}{@{} l r r r r @{}}
    \hline
       % & & \multicolumn{2}{c}{4} & \multicolumn{2}{c}{5} & \multicolumn{2}{c}{6} & \multicolumn{2}{c}{24} & \multicolumn{2}{c}{30} & \multicolumn{2}{c}{36} & \multicolumn{2}{c}{96} & \multicolumn{2}{c}{120} & \multicolumn{2}{c}{144}\\
        & & \textbf{96} & \textbf{120} & \textbf{144} \\
     %   \hline
    %    &  & STD & H & STD & H & STD & H & STD & H & STD & H & STD & H & STD & H\\
        \hline
        \multirow{2}{*}{\textbf{k-means}} & STD & -5.58 & -5.54 & -4.95 \\
                                         & R & -4.66 & -5.51 & -4.60  \\
         \multirow{2}{*}{\textbf{k-medoids}} & STD & 8.45 & 11.06 & 9.73 \\
                                            & R & 3.07 & 4.59 & 3.62 \\
         \hline
    \end{tabular}
    \label{res_simplM1_2}
\end{table}

From tables \ref{res_simplM1_1} and \ref{res_simplM1_2} it appears that the hours clustering performs the best on problem M1, especially with the range normalization and k-medoids. For approaching the reference value from below, the best approach is k-means with hours clustering. Here the range method seems slightly better than STD. 

k-medoids constantly overestimates the objective value while k-means constantly underestimates it. However, in general and for around 30 days and 120 hours, the k-means seems to be the appropriate choice. Indeed, even though k-medoids with STD also has good results, it appears less consistent. With this algorithm, the performance does not always improve with an increasing number of clusters; choosing the correct amount of clusters would become harder. k-means, while not completely exempt from this flaw, appears more robust in this regard.

For M1, the average of the run time for days clustering for 24, 30 and 36 days is 3500 seconds with extreme values of 2 289 and 5628 seconds. For the hours clustering, the average runtime is 5 973 seconds with extremes of 2 421 and 12 500 seconds. Days clustering is on average almost twice as fast as hours clustering on this simplified model despite having more timesteps overall. As a reference, to solve the problem without any clustering (using a complete year) takes around 30 000 seconds.

For M0 the runtimes are low with all values below 360 seconds.

\subsection{Complete Model}

For the complete model, no reference value is presented because running the models with a complete year of data takes too long and it is the reason clustering is explored in the first place.

\begin{figure}[h]
    \centering
    \includegraphics[width=0.48\textwidth,trim=4 4 4 4,clip]{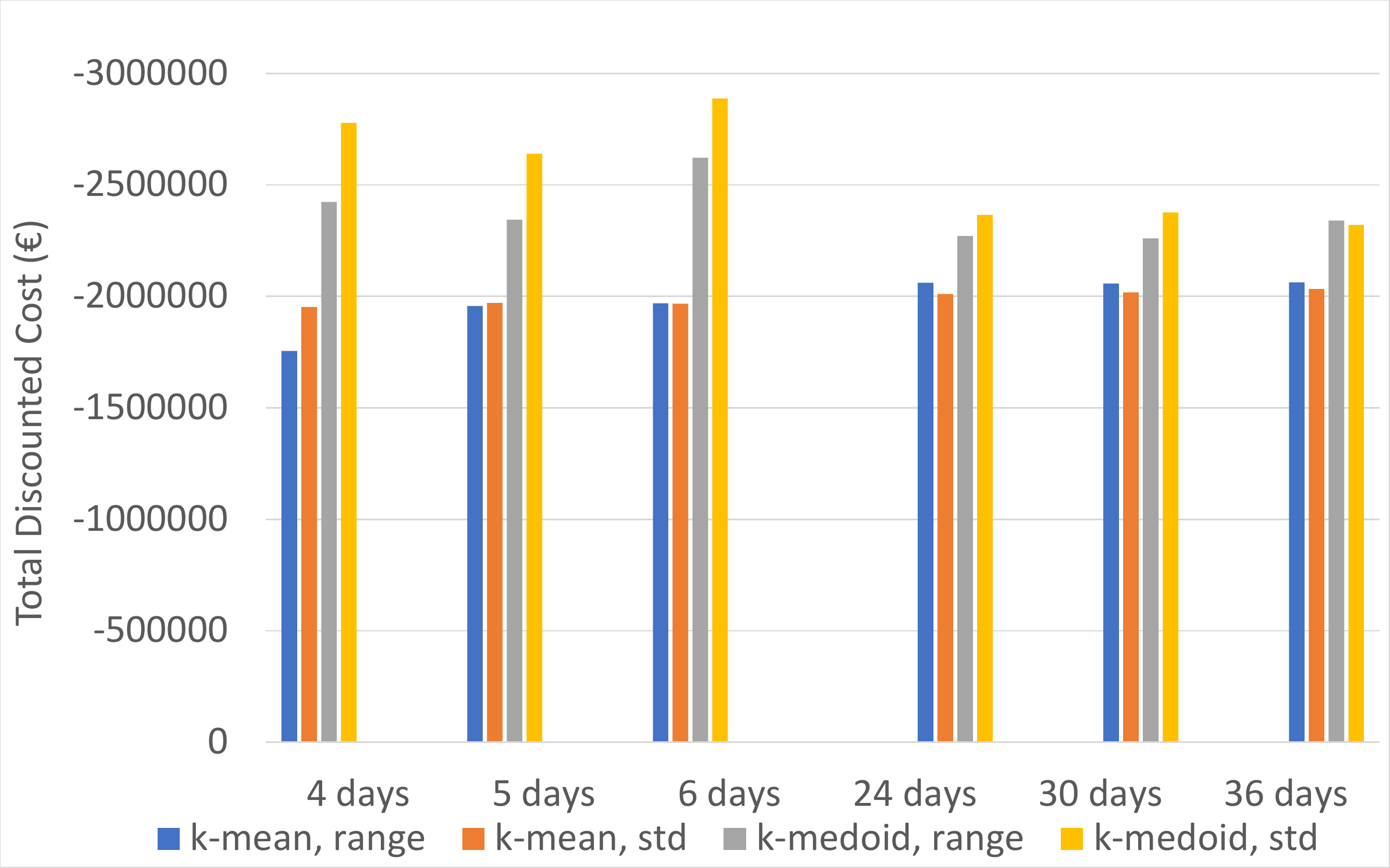}
    \caption{Objective Values for M0 with Design Days with Complete Model}
    \label{fig:objval_fullmodel}
\end{figure}

 Figure \ref{fig:objval_fullmodel} presents the objective values resulting from the optimization in the case of M0. Without a reference value, it is impossible to reach a conclusion regarding the performances of each approach. However we can make some remarks. The objective values follow the same patterns as in the case of the simplified model and from the results we can expect that in this case as well k-means underestimates and k-medoids overestimates the objective value. It also appears that even a few days are enough to get satisfying results when using k-means.
 
\begin{table}[t]
\renewcommand{\arraystretch}{1.15}
    \centering
    \caption{Runtime for M0 in seconds with days (\textbf{STD}: Standard Devation, \textbf{R}:Range)}
    \begin{tabular}{@{} l r r r r r r r  @{}}
    \hline
       % & & \multicolumn{2}{c}{4} & \multicolumn{2}{c}{5} & \multicolumn{2}{c}{6} & \multicolumn{2}{c}{24} & \multicolumn{2}{c}{30} & \multicolumn{2}{c}{36} & \multicolumn{2}{c}{96} & \multicolumn{2}{c}{120} & \multicolumn{2}{c}{144}\\
         & &  \textbf{4} & \textbf{5} & \textbf{6} & \textbf{24} & \textbf{30} & \textbf{36}\\
     %   \hline
    %    &  & STD & H & STD & H & STD & H & STD & H & STD & H & STD & H & STD & H\\
        \hline
        \multirow{2}{*}{\textbf{k-means}}  & STD & 70.98 & 96.47 & 108.8 & 1708 & 2846 & 3342  \\
                                          & R & 116.3 & 115.2 & 155.4 & 1924 & 3504 & 4320  \\
        \multirow{2}{*}{\textbf{k-medoids}} & STD & 63.42 & 62.98 & 288.4 & 3544 & 4157 & 6163 \\
                                           & R & 57.52 & 115.0 & 127.5 & 2088 & 3442 & 5288  \\
        \hline
    \end{tabular}
    \label{res_M0}
\end{table}

Regarding runtime for M0, k-means with STD is clearly the fastest while k-medoids with STD is the slowest being about half as fast. k-medoids range and k-means range have comparable runtimes except for the case of 36 days where the k-medoids version is about 20\% slower. k-means range is itself 25\% slower than k-means STD.

\begin{figure}[h]
    \centering
    \includegraphics[width=0.48\textwidth,trim=4 4 4 4,clip]{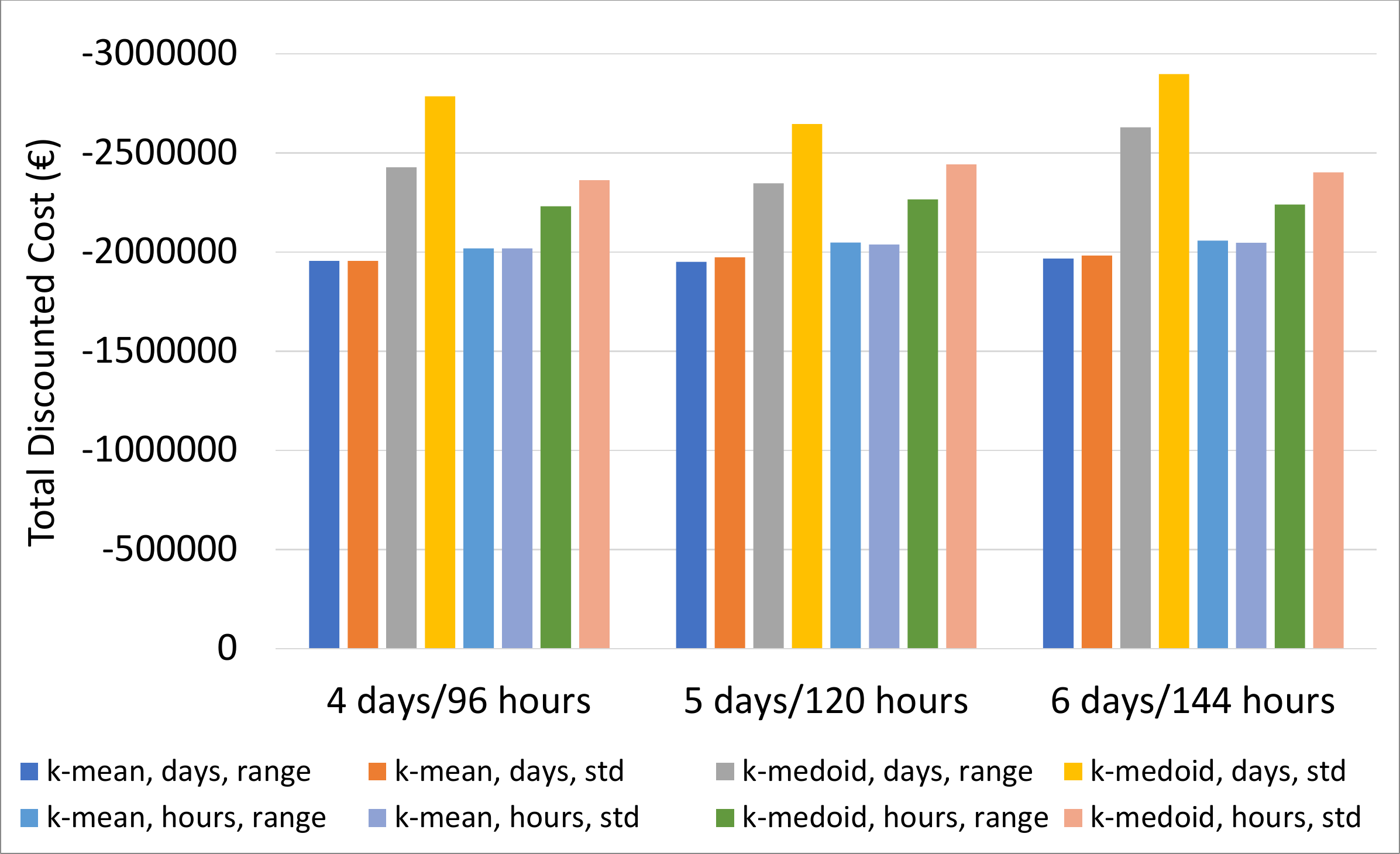}
    \caption{Objective Values for M1 with Days and Hours with Complete Model}
    \label{fig:objvalM1_fullmodel}
\end{figure}

For M1, the same remarks hold. K-medoids and k-means seem to respectively over- and underestimate the objective value. Figure \ref{fig:objvalM1_fullmodel} confirms that for k-medoids, the range method performs better than STD as in Table \ref{res_simplM1_1} and \ref{res_simplM1_2}. For k-means we also find that the results are similar.

\begin{table}[t]
\renewcommand{\arraystretch}{1.15}
    \centering
    \caption{Runtime for M1 in seconds(\textbf{STD}: Standard Devation, \textbf{R}:Range)}
    \begin{tabular}{@{} l r r r r @{}}
    \hline
       % & & \multicolumn{2}{c}{4} & \multicolumn{2}{c}{5} & \multicolumn{2}{c}{6} & \multicolumn{2}{c}{24} & \multicolumn{2}{c}{30} & \multicolumn{2}{c}{36} & \multicolumn{2}{c}{96} & \multicolumn{2}{c}{120} & \multicolumn{2}{c}{144}\\
        & \textbf{Days} &  \textbf{4} & \textbf{5} & \textbf{6} \\
     %   \hline
    %    &  & STD & H & STD & H & STD & H & STD & H & STD & H & STD & H & STD & H\\
        \hline
        \multirow{2}{*}{\textbf{k-means}} & STD & 1386 & 2250 & 4340   \\
                                         & R & 2059 & 3393 & 3519    \\
         \multirow{2}{*}{\textbf{k-medoids}} & STD & 1723 & 1717 & 5838 \\
                                            & R & 1139 & 2319 & 2626  \\
        \hline
        & \textbf{Hours} & \textbf{96} & \textbf{120} & \textbf{144}\\
        \hline
        \multirow{2}{*}{\textbf{k-means}} & STD & 14789 & 29159 & 62239  \\
                                         & R & 13342 & 45048 & 60165   \\
         \multirow{2}{*}{\textbf{k-medoids}} & STD & 20288 & 55509 & 105672  \\
                                            & R & 18632 & 19860 & 59055 \\
         \hline
    \end{tabular}
    \label{res_M1_1}
\end{table}

M1 is between 15 and 40 times longer to solve than M0 for the days. When it comes to the difference between the days and the hours, even though the number of timesteps are the same, the hourly model takes at least 10 times longer to solve than the daily model. This difference is hard to explain. Indeed both models get the same number of timesteps and are identical with the exception of what is presented in equations \ref{hour} and \ref{day}.

If the use of k-medoids is required for any reason, then using the hourly method can bring significant improvements to the precision over the daily method. These improvements needs to be considered in regard to the increased solving time to choose the method to use. Otherwise, k-means should be preferred. In that case, the improvements of the precision is insufficient to justify using the hourly method.
One such possible reason is to have a good representation of the solar irradiance which is the case for ZENIT. By using the day method with low numbers of days, even though the solving time and objective values are good, the representation of the solar irradiance is problematic as seen in Fig. \ref{irr}. In our case and to get a good solar irradiance representation, the use of k-means and hours clustering in M1 is preferable.

%While the use of M0 or M1 should be considered on the basis of the necessity to include seasonal storage, the choice of the clustering method (algorithm, cluster type and normalization method) can be made based on the results presented in this paper. For the particular application of designing the energy system of neighborhoods with an objective of zero emission, the best methods appears to be to use the k-means algorithm with the range normalization and days as cluster type. A low number of days is fine but it can be interesting to increase. The trade-off between time and precision should then be considered.
%\section{Assessment of clustering for the purpose of the MPC}
%
%Maybe should be part of separate paper along with the different strategies for ZENIT %in MPC

\section{Limitations}

There are different limitations that should be mentioned regarding this paper. Regarding the studied methods, the fact that only clustering algorithms are studied have been explained; however other clustering algorithms could offer advantages. %Another limitation is that only one clustering is explored.
Many heuristics, either new or variations around the one used, could also be studied and finding the overall best heuristic presents a challenge. 
The clustering has been used on a specific case and we cannot guarantee that the same result holds true for larger cases or in other countries where the correlation between the different inputs are different. 
Unfortunately no reference value is shown for the complete model and a simplified model had to be used in order to compare the precision.

\section{Conclusion}

After introducing the use of reduction techniques and clustering in energy systems and in particular in the design of the energy system of neighborhoods, this paper discussed why clustering is chosen over other solutions such as downsampling. Different clustering methods have then been evaluated, first directly on their ability to come close to the original dataset and then on the results they give when used in ZENIT. 
K-means and k-medoids have been compared and the study allowed to highlight that counter to what is found for many other energy system applications, k-means performs better than k-medoids. 
The study also highlights the role of the normalization method on the performances by comparing a method using the standard deviation and one using the range of values. 
We find occurrences of models using clustered days (or design days) and of instances using clustered hours in the literature but the reason for the choice are not always clear. In this study, both approaches are implemented and the relation between the performance, the solving time and the possible uses of each are reviewed.
The impact of the use of a simple heuristic is also studied.
Two versions of the optimization models were used with different capabilities when it comes to storage: M0 for daily storage operation and M1 for storage without time limitation.
While the use of M0 or M1 should be considered on the basis of the necessity to include seasonal storage, the choice of the clustering method (algorithm, cluster type and normalization method) can be made based on the results presented in this paper. For the particular application of designing the energy system of neighborhoods with an objective of zero emissions, the best method appears to be to use the k-means algorithm with the range normalization and days as cluster type. A low number of days is fine but it can be interesting to increase it to improve the representation of the solar irradiance for example. The trade-off between time and precision should then be considered.
Further work could extend the result to other cases and study if the results presented in this paper scale to bigger neighborhoods. Other clustering algorithms or heuristics could also be investigated.

\section*{Acknowledgment}
This article has been written within the Research Center on Zero Emission Neighborhoods in Smart Cities (FME ZEN). The author gratefully acknowledges the support from the ZEN partners and the Research Council of Norway.

The author would also like to thank John Clau{\ss} for providing the hourly $CO_2$ factor for electricity data.

\bibliographystyle{IEEEtran}
\bibliography{biblio}

% or


\appendices

\section{Additional Results of the Clustering}

Additional results are presented in this appendix. In particular, the RMSD for the time series that were not included in section \ref{clust} are shown in this section.

\begin{figure}[h]
    \centering
    \subfloat[Hours]{\includegraphics[width=0.24\textwidth]{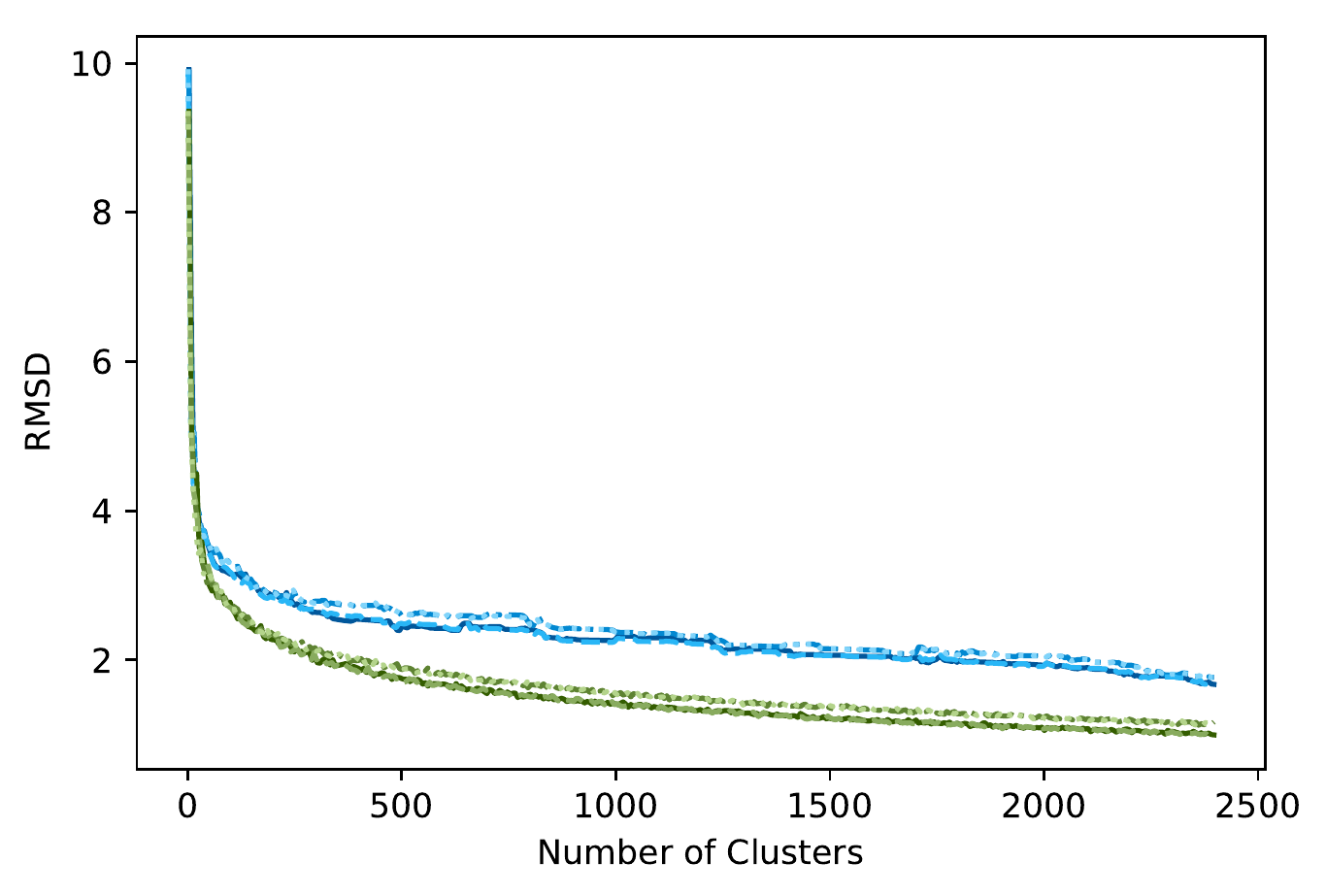}}
    \subfloat[Days]{\includegraphics[width=0.24\textwidth]{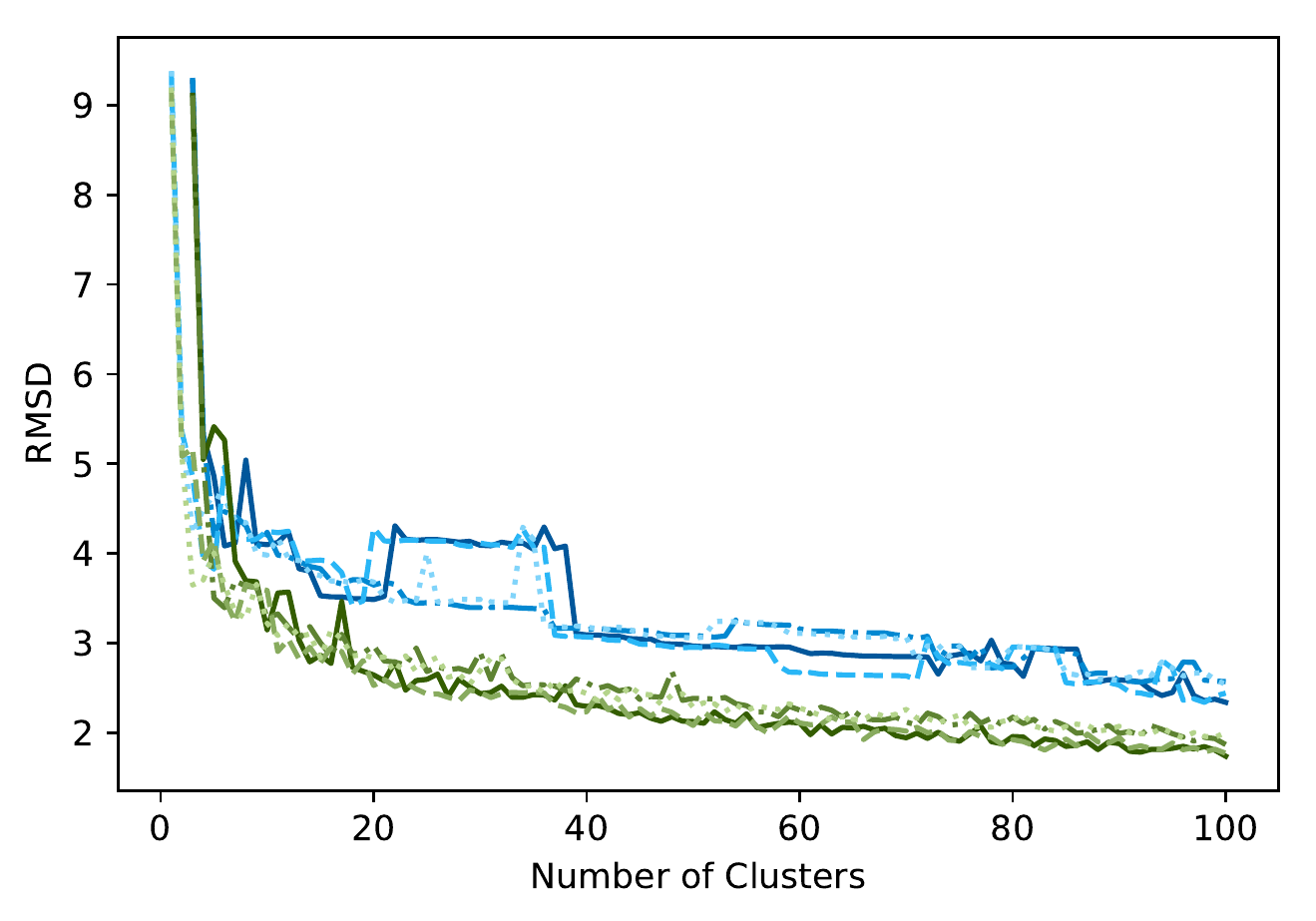}}
    \caption{RMSD of Temperature}
    \label{temperature}
\end{figure}

The errors for the temperature time series, Fig. \ref{temperature}, are very similar to the overall ones. The RMSD of temperature plateaus rather quickly to around 2 for the hours, and 2.8 for the days.

\begin{figure}[h]
    \centering
    \subfloat[Hours]{\includegraphics[width=0.24\textwidth]{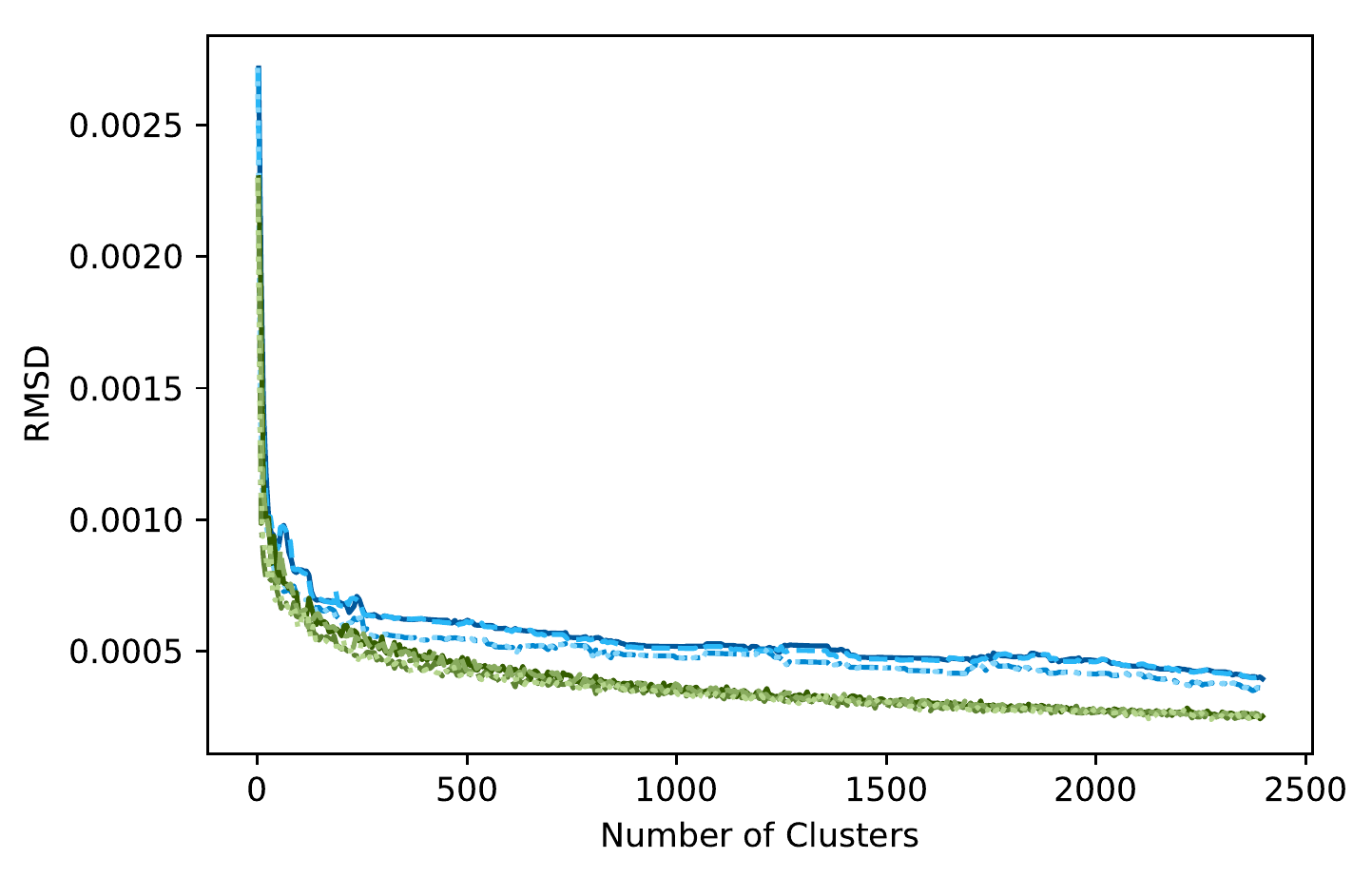}}
    \subfloat[Days]{\includegraphics[width=0.24\textwidth]{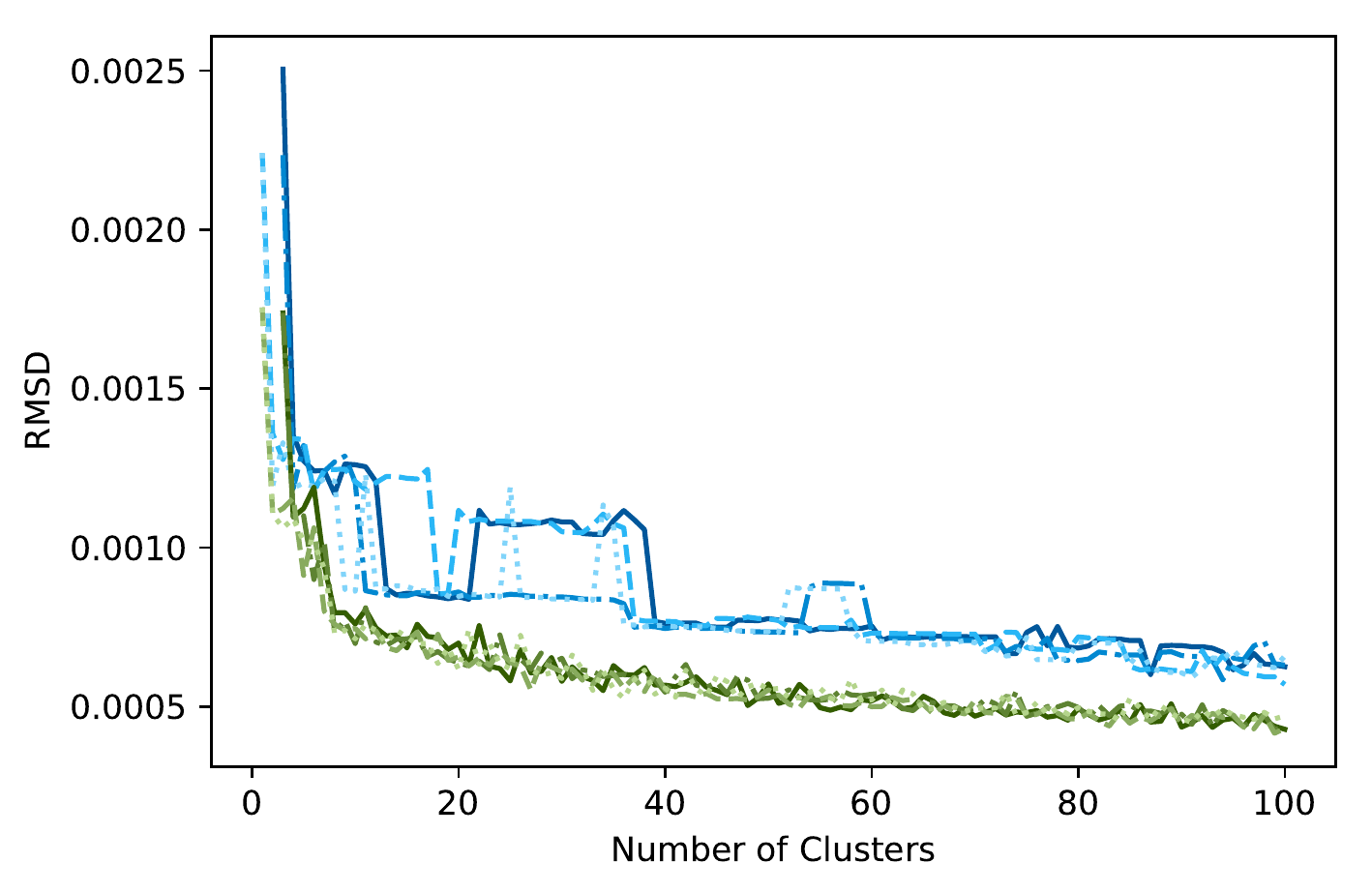}}
    \caption{RMSD of DHW Load in the Normal Offices}
\end{figure}

\begin{figure}[h]
    \centering
    \subfloat[Hours]{\includegraphics[width=0.24\textwidth]{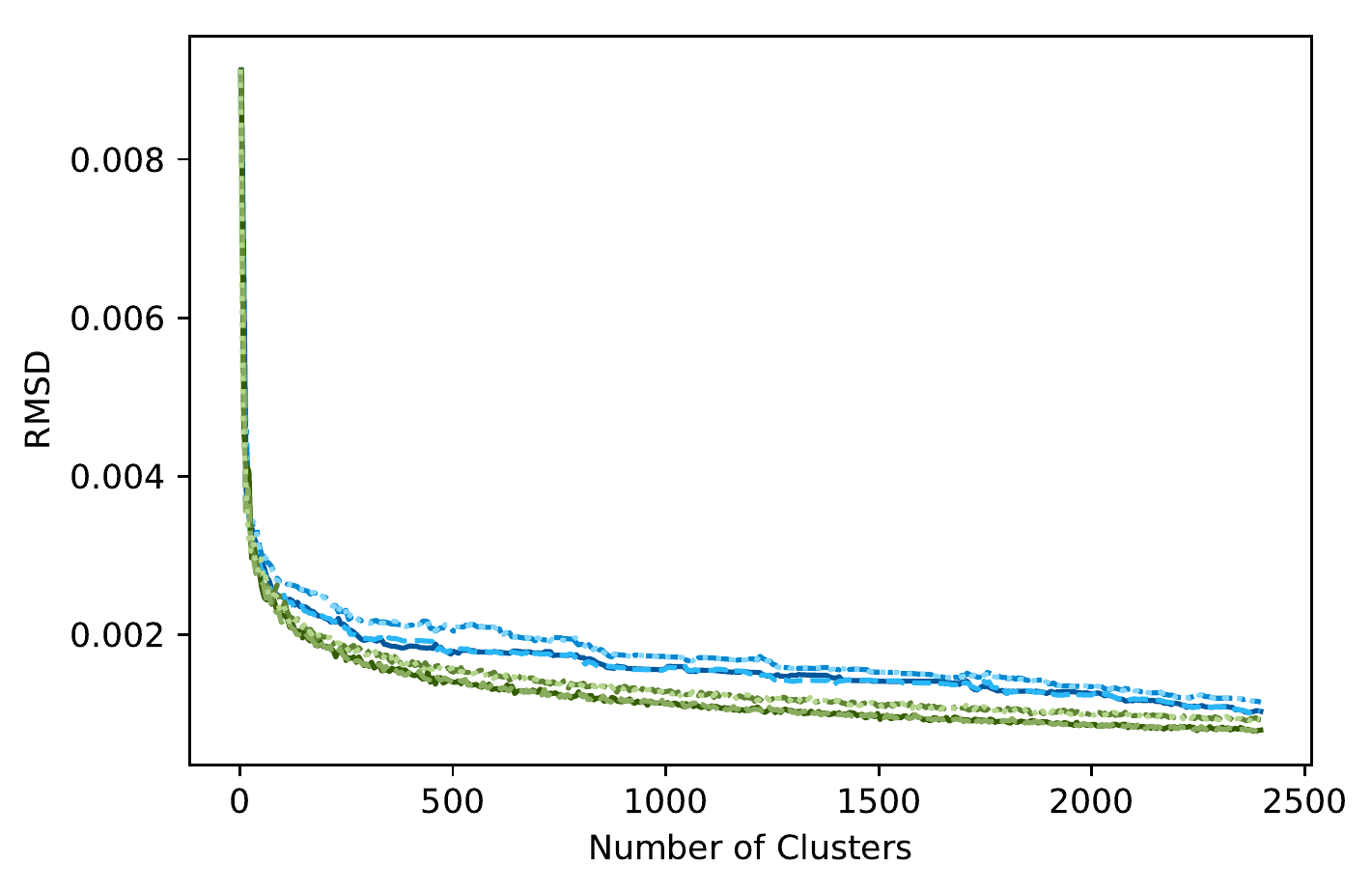}}
    \subfloat[Days]{\includegraphics[width=0.24\textwidth]{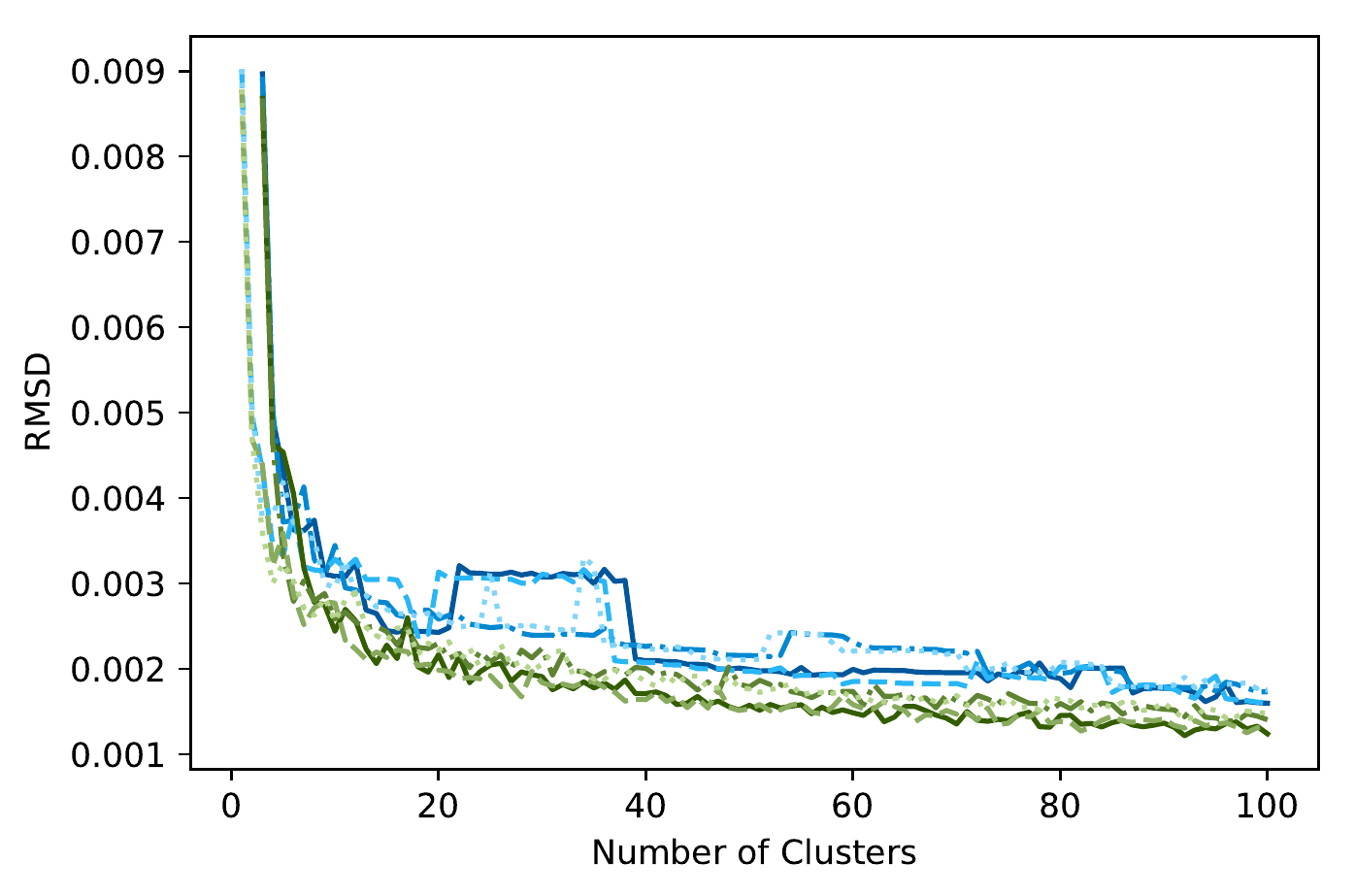}}
    \caption{RMSD of SH Load in the Normal Offices}
\end{figure}

For the offices already at the passivhus standard:
 
\begin{figure}[h]
    \centering
    \subfloat[Hours]{\includegraphics[width=0.24\textwidth]{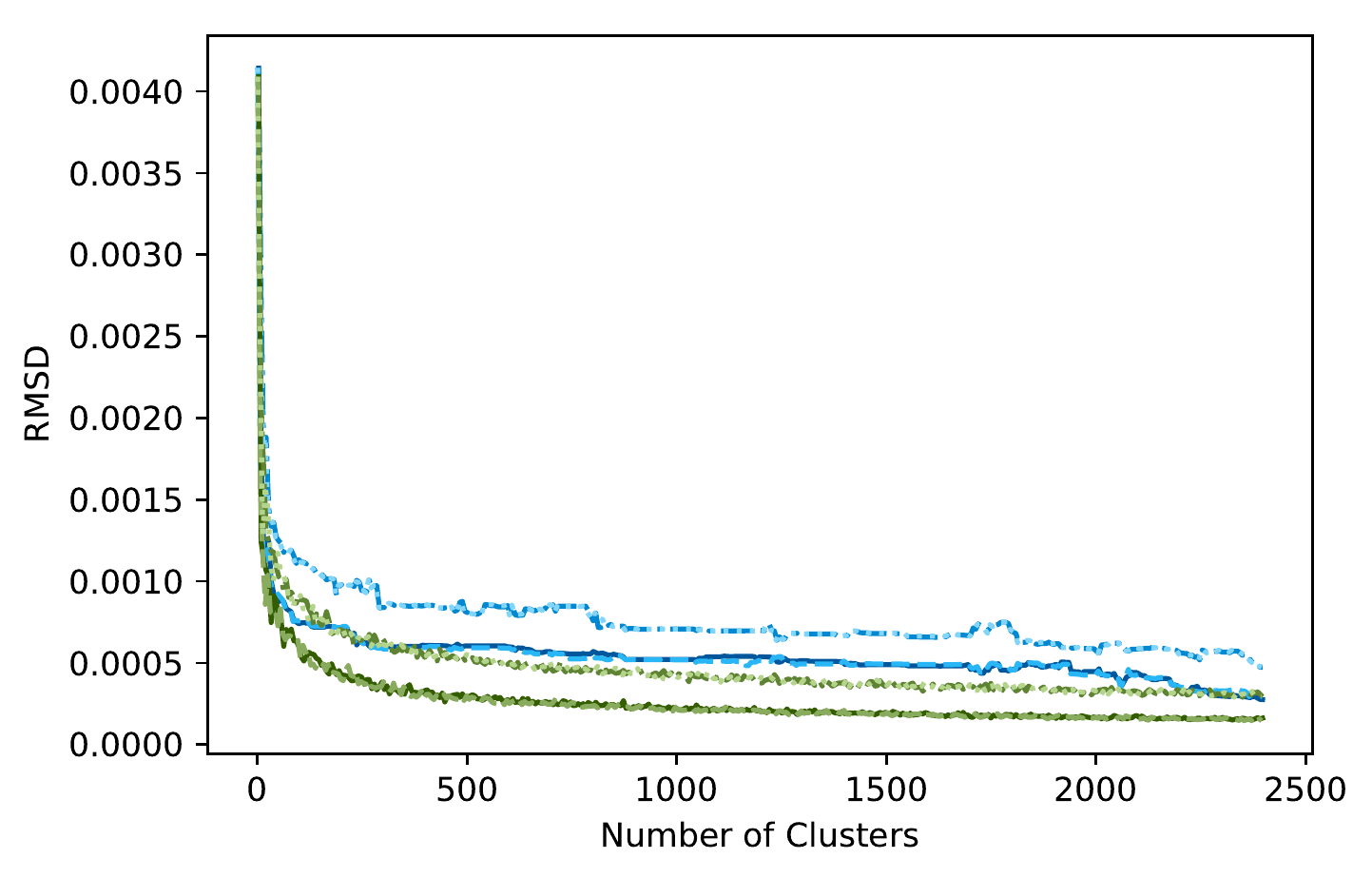}}
    \subfloat[Days]{\includegraphics[width=0.24\textwidth]{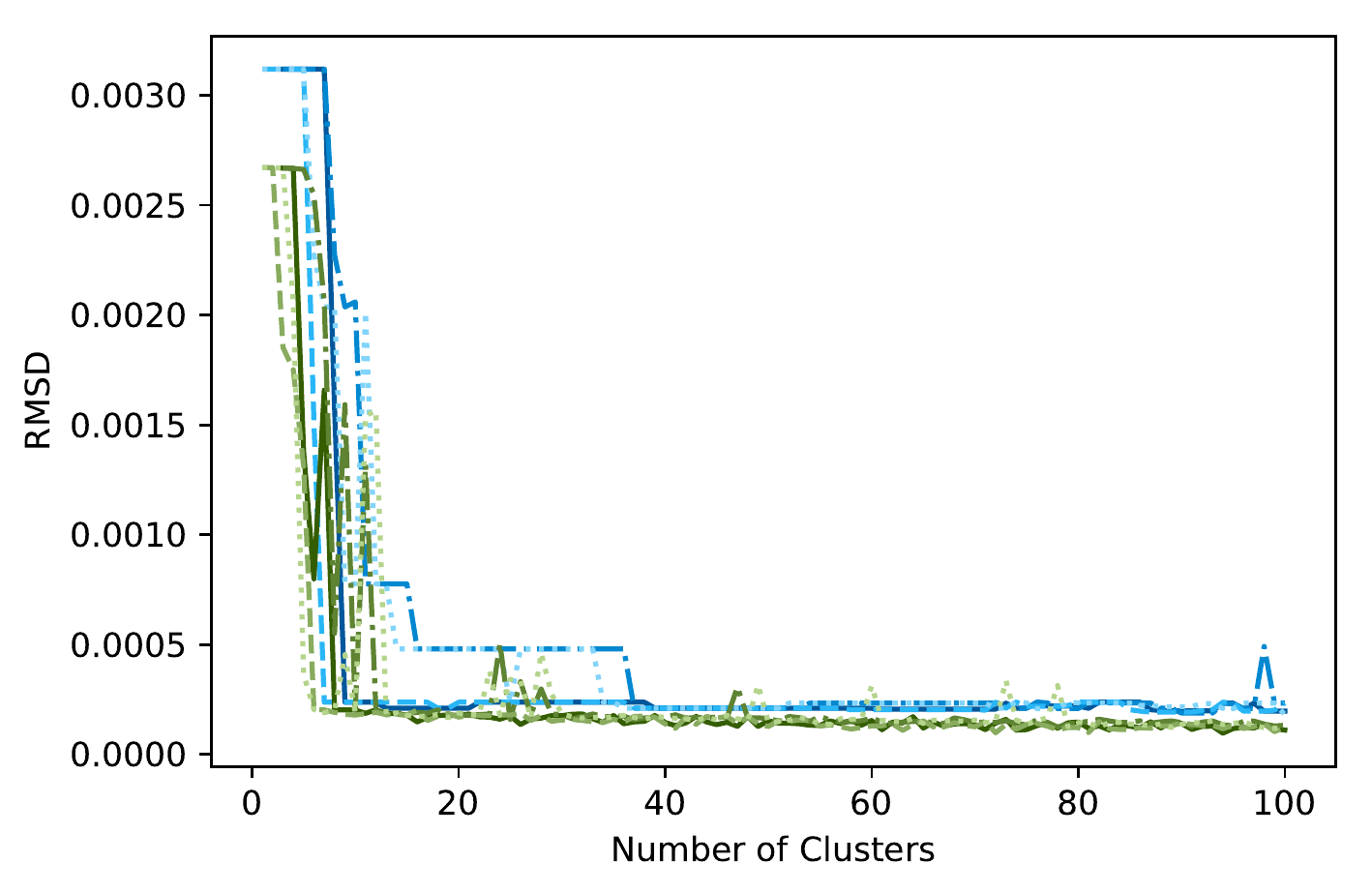}}
    \caption{RMSD of Electric Load in the Passive Offices}
\end{figure}

\begin{figure}[h]
    \centering
    \subfloat[Hours]{\includegraphics[width=0.24\textwidth]{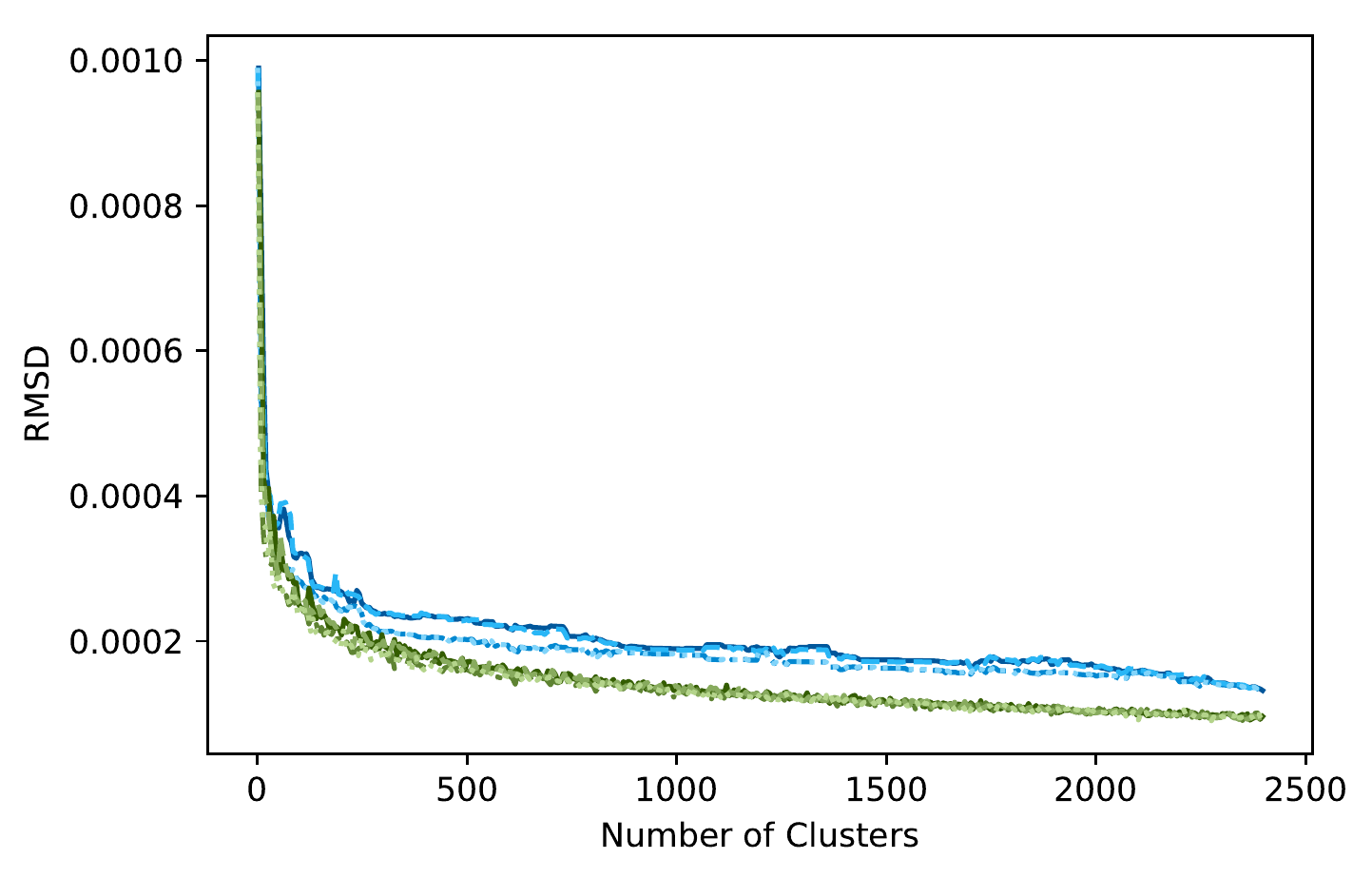}}
    \subfloat[Days]{\includegraphics[width=0.24\textwidth]{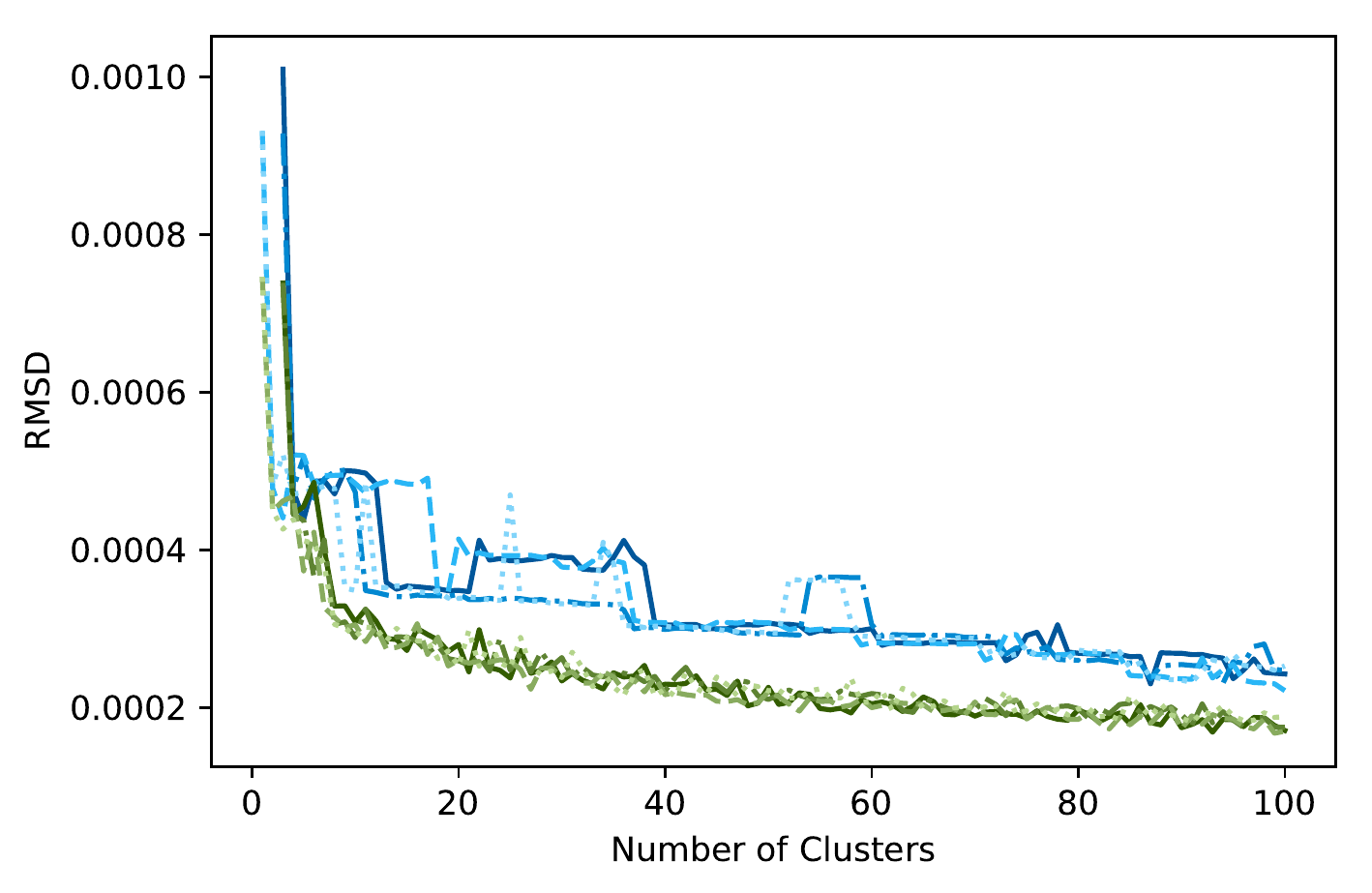}}
    \caption{RMSD of DHW Load in the Passive Offices}
\end{figure}

\begin{figure}[h]
    \centering
    \subfloat[Hours]{\includegraphics[width=0.24\textwidth]{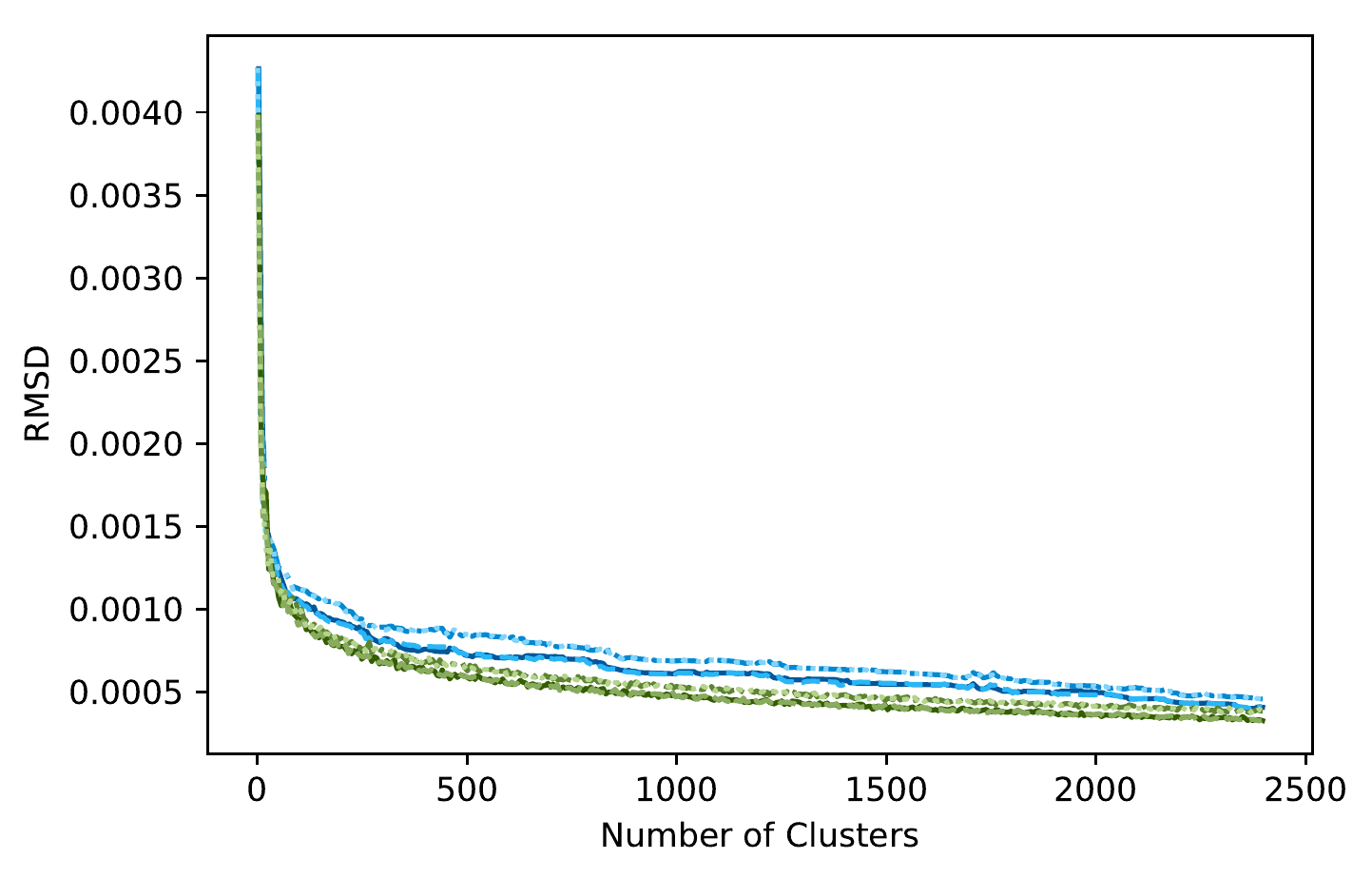}}
    \subfloat[Days]{\includegraphics[width=0.24\textwidth]{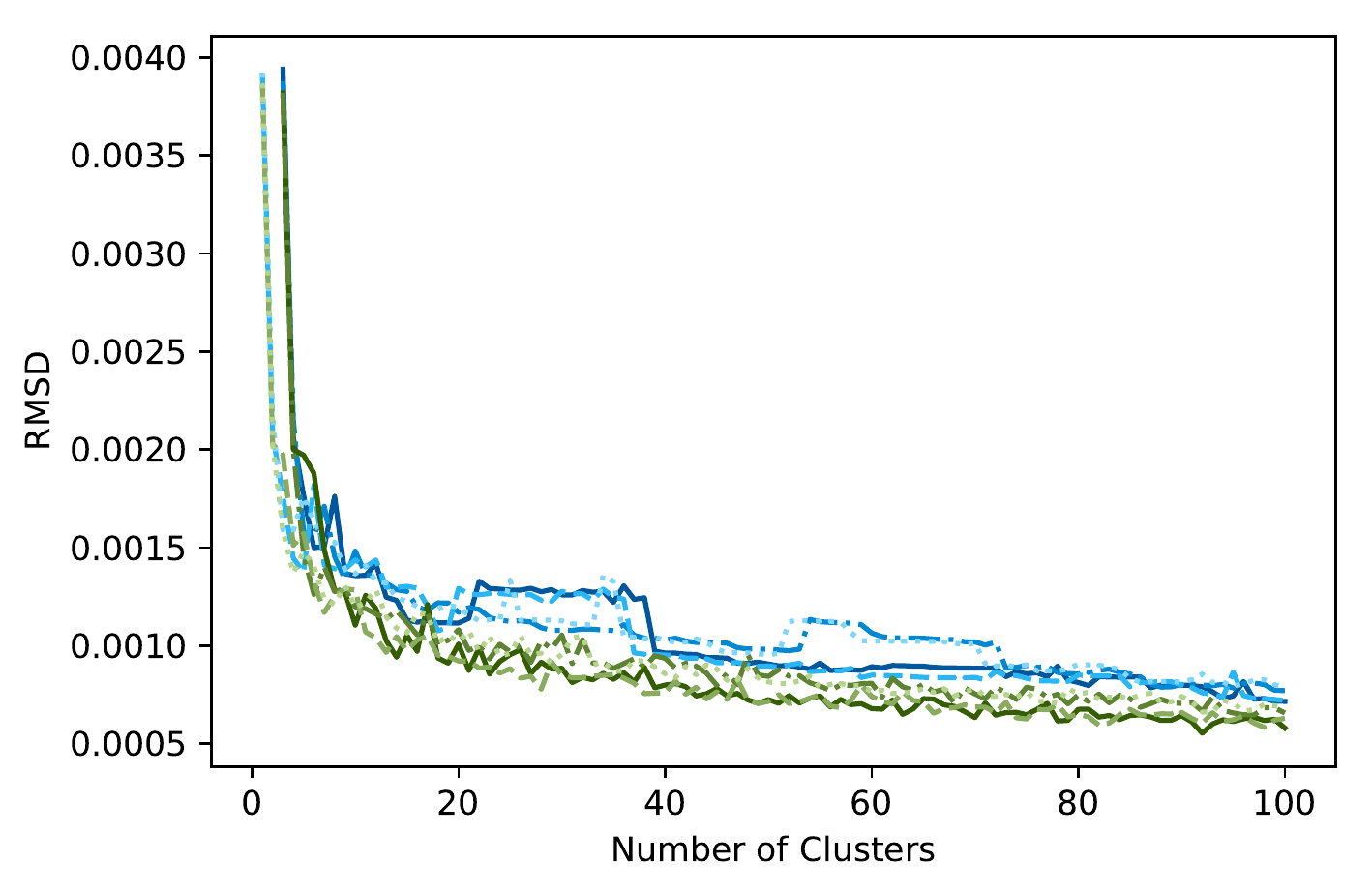}}
    \caption{RMSD of SH Load in the Passive Offices}
\end{figure}

For the student housings:

\begin{figure}[h]
    \centering
    \subfloat[Hours]{\includegraphics[width=0.24\textwidth]{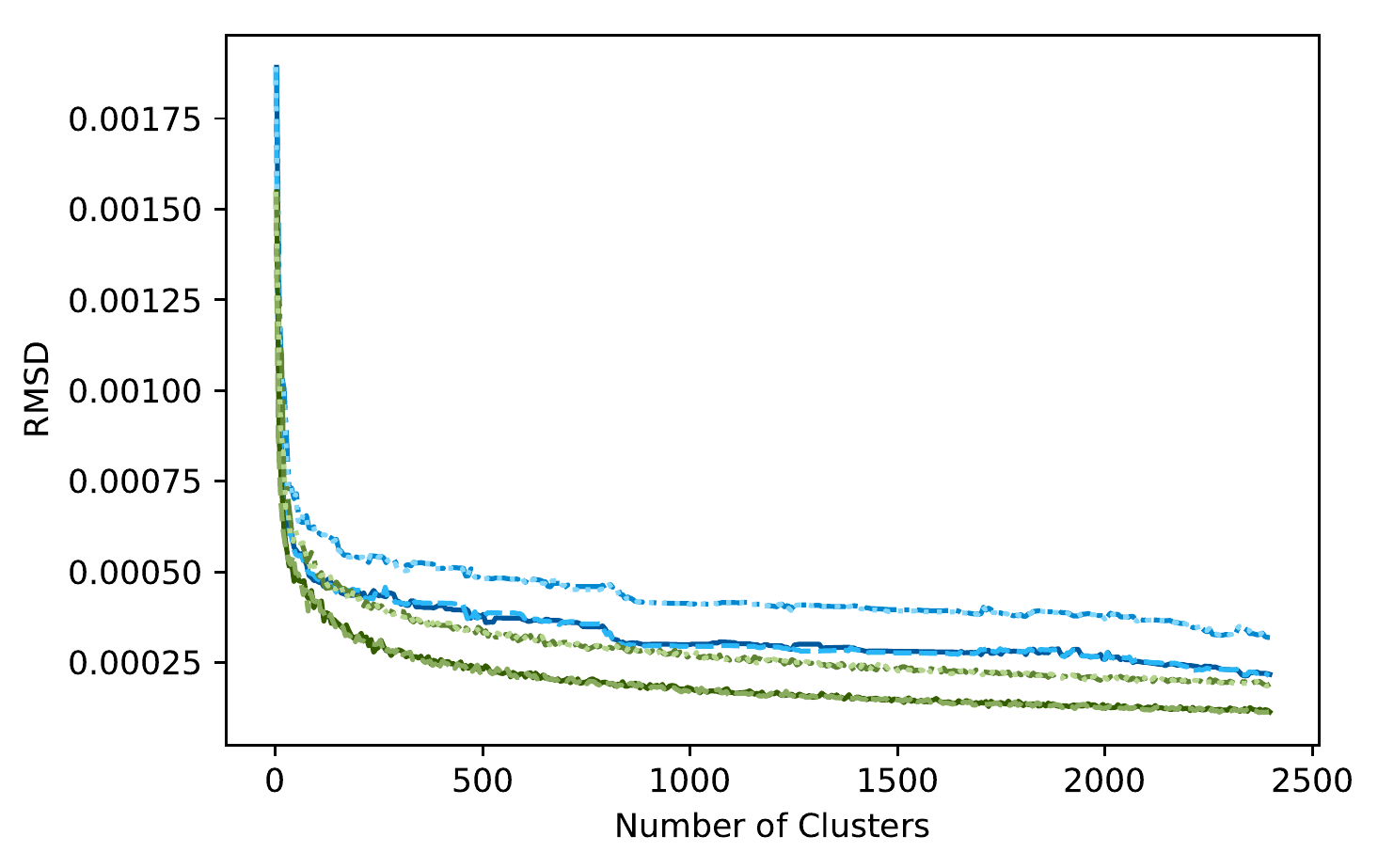}}
    \subfloat[Days]{\includegraphics[width=0.24\textwidth]{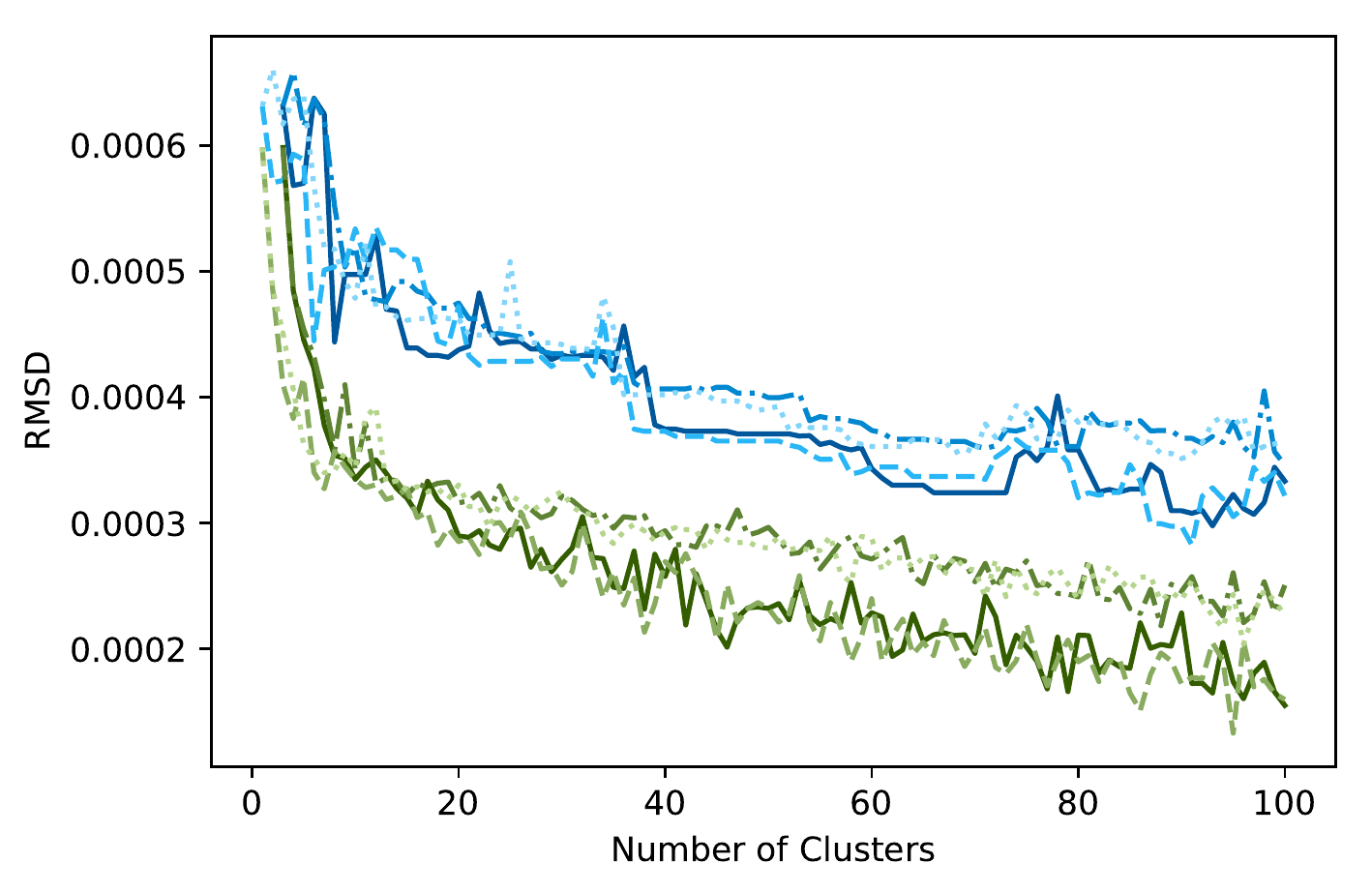}}
    \caption{RMSD of Electric Load in the Student Housing}
\end{figure}

\begin{figure}[h]
    \centering
    \subfloat[Hours]{\includegraphics[width=0.24\textwidth]{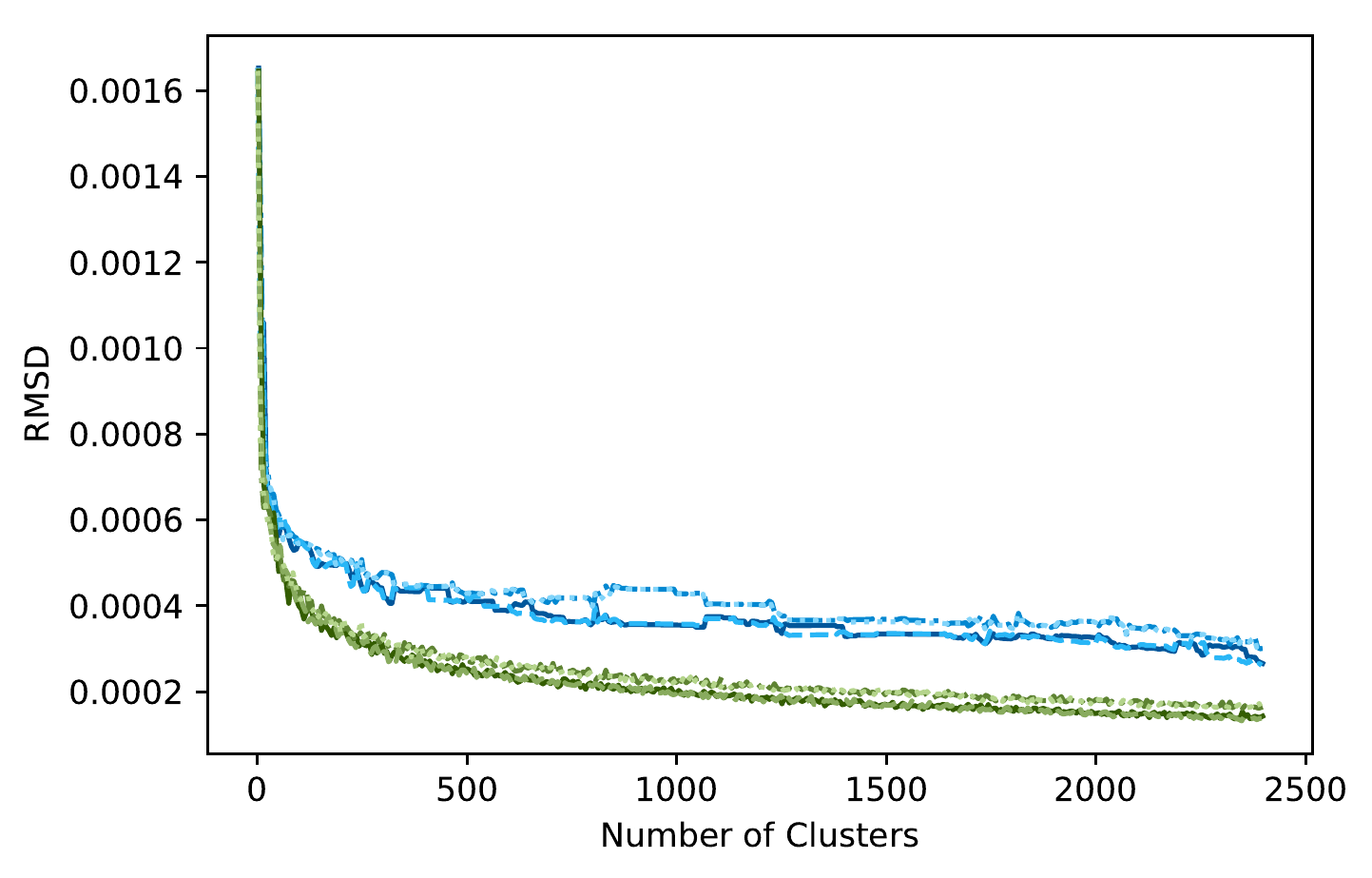}}
    \subfloat[Days]{\includegraphics[width=0.24\textwidth]{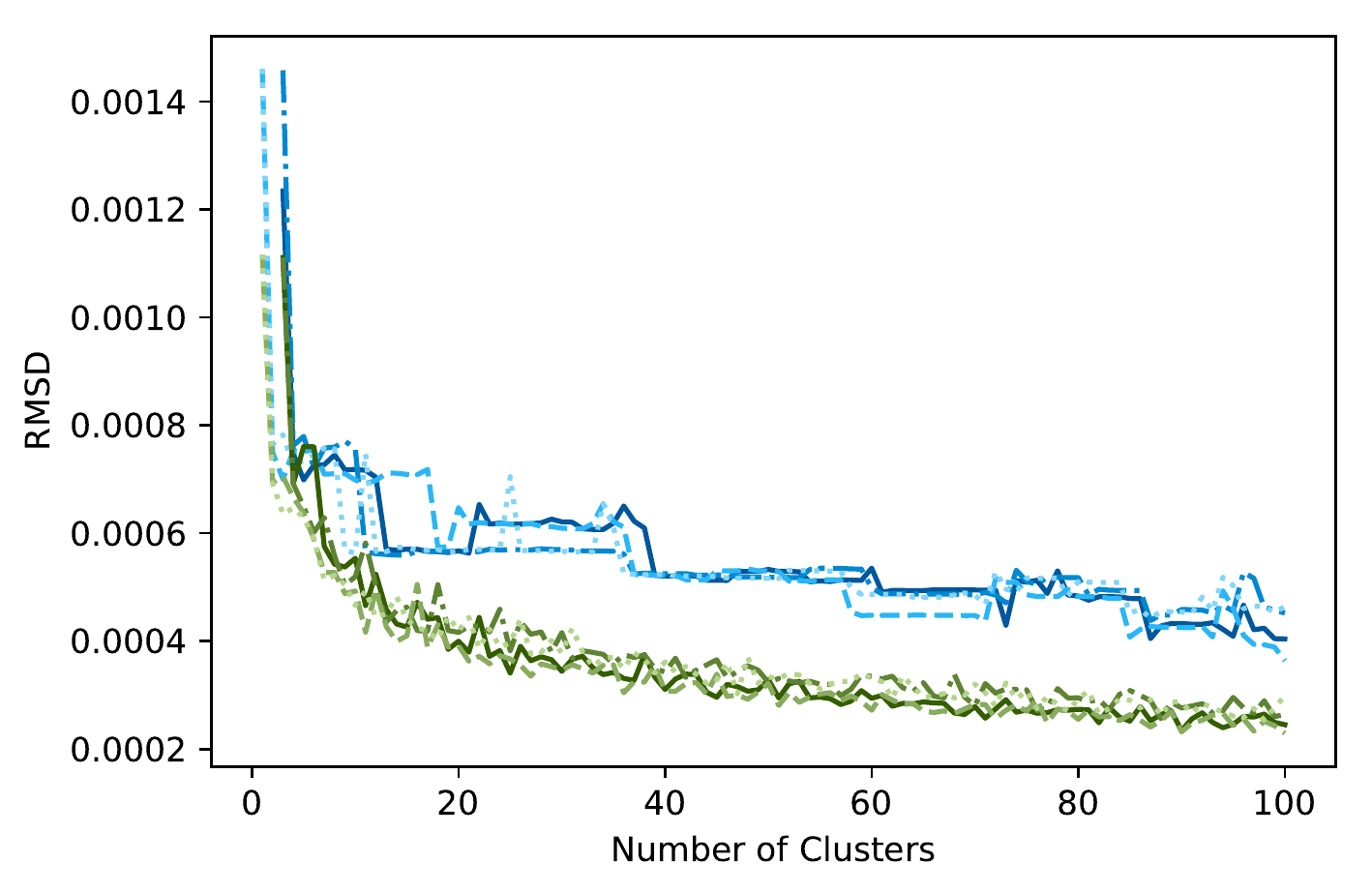}}
    \caption{RMSD of DHW Load in the Student Housing}
\end{figure}

\begin{figure}[h]
    \centering
    \subfloat[Hours]{\includegraphics[width=0.24\textwidth]{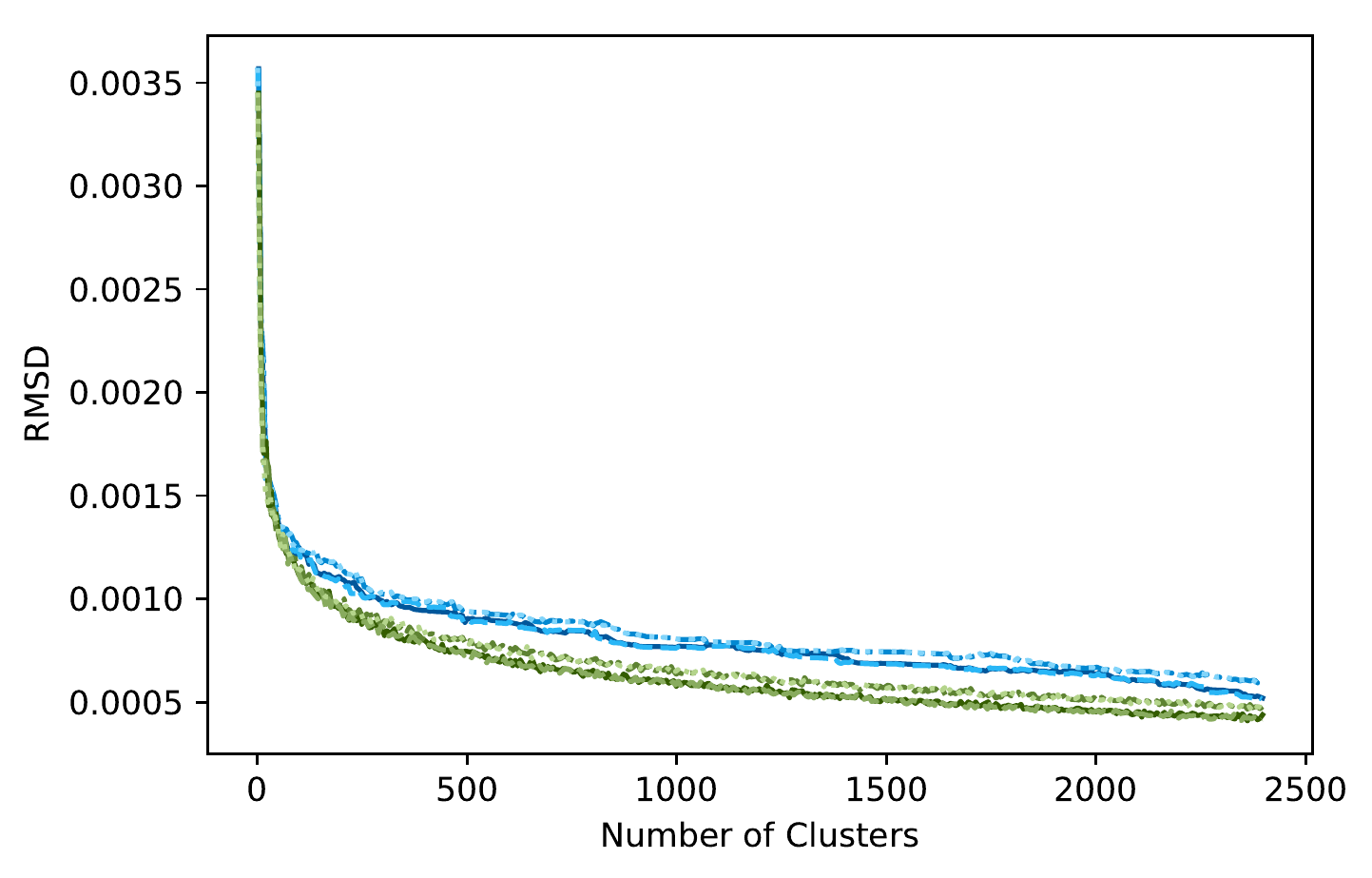}}
    \subfloat[Days]{\includegraphics[width=0.24\textwidth]{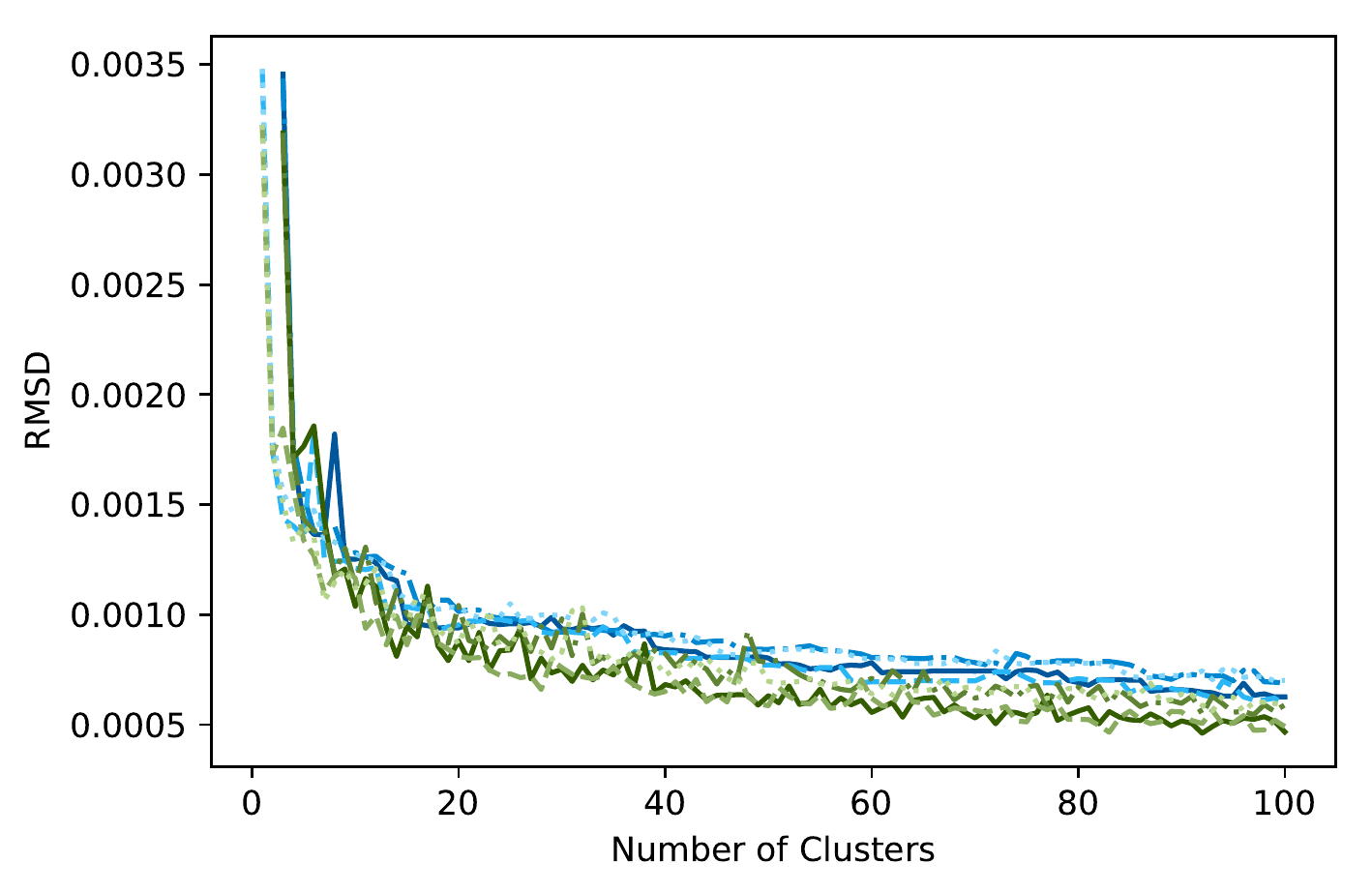}}
    \caption{RMSD of SH Load in the Student Housing}
\end{figure}

The figures for the yearly average errors presented in Table \ref{tab:cluster_res} are presented here as well.

\begin{figure}[h]
    \centering
    \subfloat[Hours]{\includegraphics[width=0.24\textwidth]{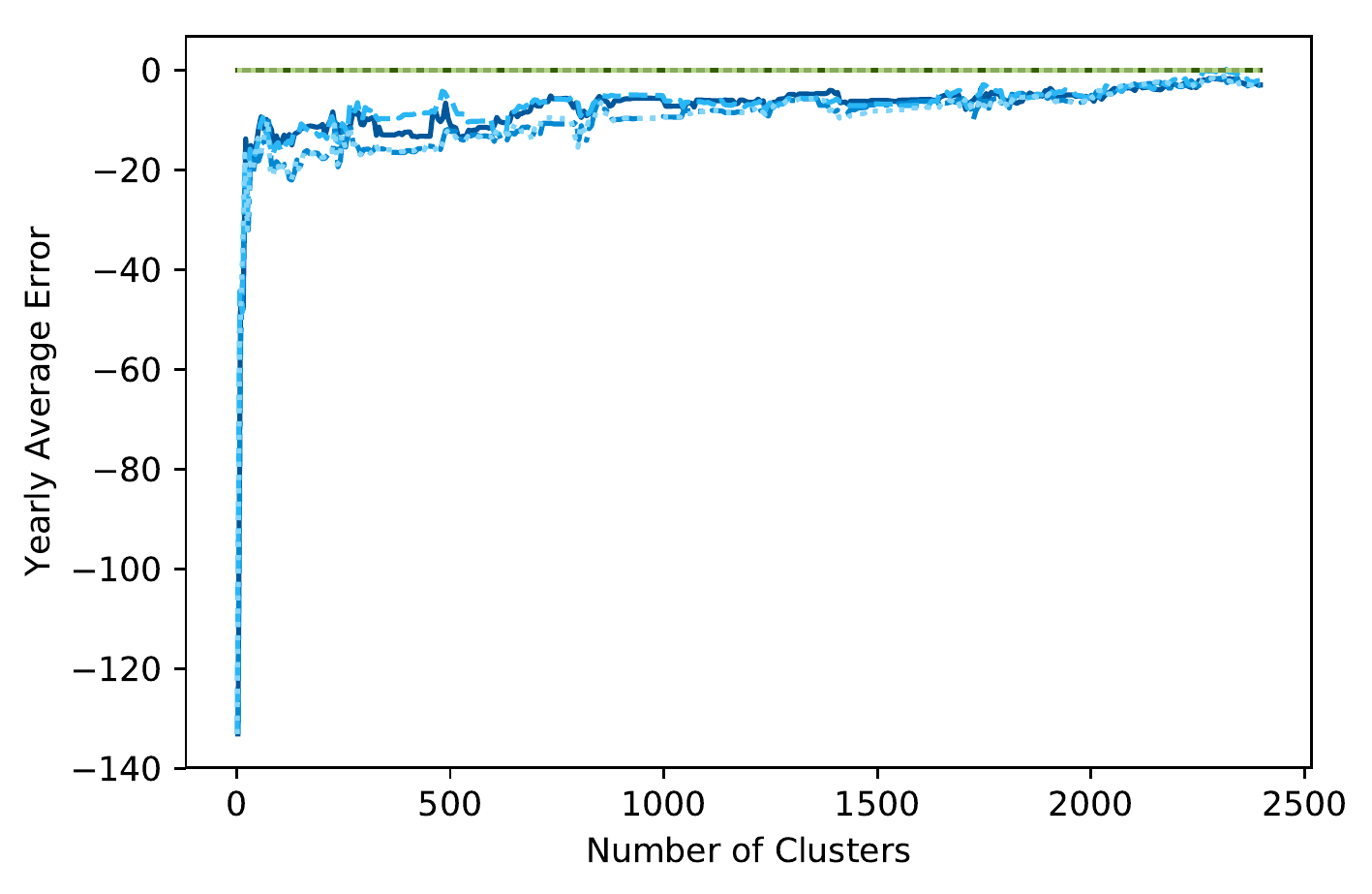}}
    \subfloat[Days]{\includegraphics[width=0.24\textwidth]{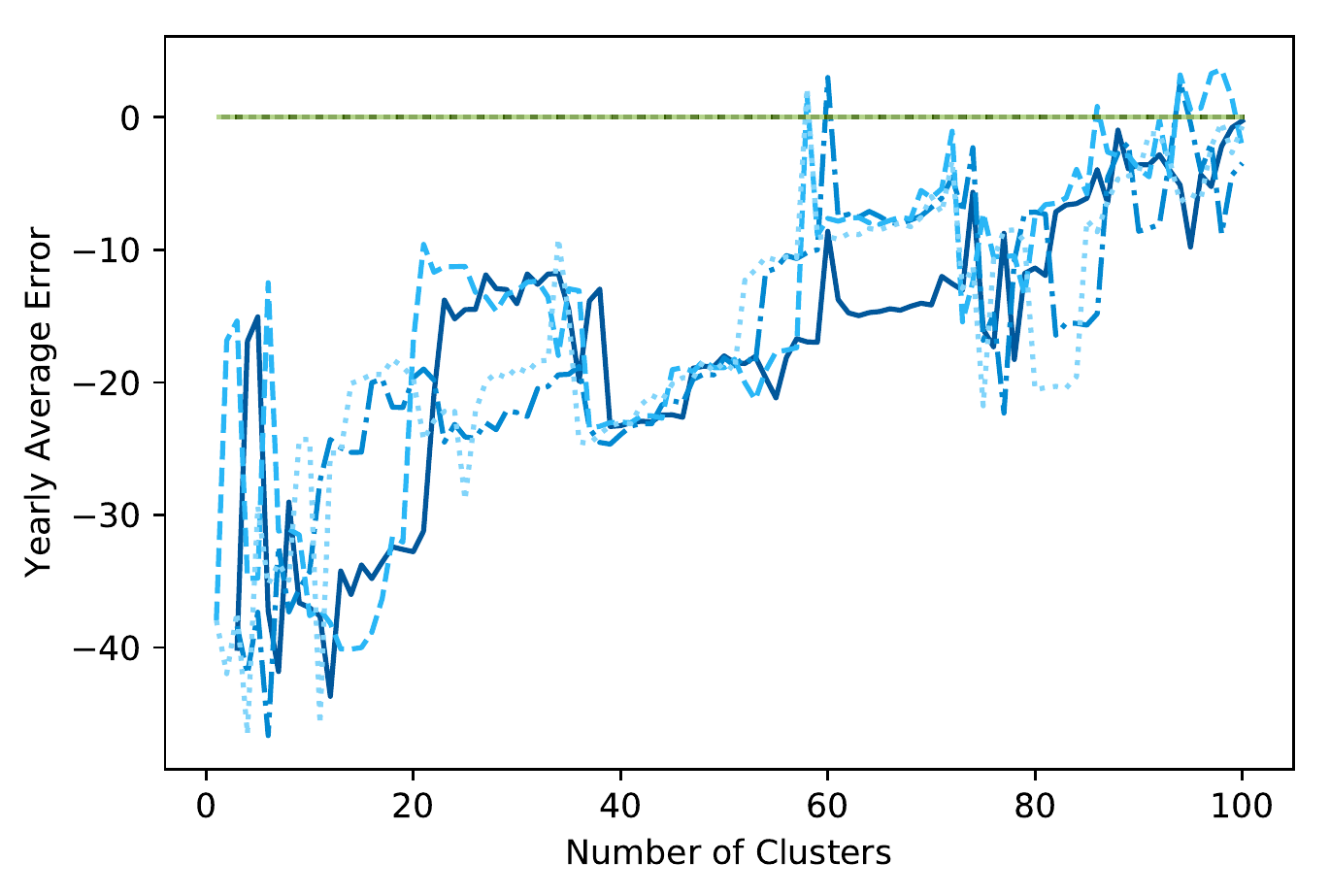}}
    \caption{YAE of Irradiance}
\end{figure}

%\section*{Acknowledgment}
%This article has been written within the Research Center on Zero Emission %Neighborhoods in Smart Cities (FME ZEN). The authors gratefully acknowledge the %support from the ZEN partners and the Research Council of Norway.
%
%The authors would also like to thank John Clau{\ss} for providing the hourly $CO_2$ %factor for electricity data.

\ifCLASSOPTIONcaptionsoff
  \newpage
\fi

\end{document}